\documentclass[preprint,twoside,11pt]{article}

\usepackage{jmlr2e_number}
\usepackage{blindtext}
\usepackage{additional_package}
\usepackage{amsmath,amsfonts,amssymb}
\usepackage[T1]{fontenc}
\usepackage{xcolor} 
\usepackage{mathrsfs} 
\usepackage{quiver} 
\usepackage{enumitem}
\usepackage{hyperref}
\usepackage{thm-restate} 
\usepackage{subcaption}

\usepackage{lastpage}
\jmlrheading{25}{2025}{1-\pageref{LastPage}}{}{}{}{Salem Said, Franziskus Steinert, and Cyrus Mostajeran}
\ShortHeadings{title}{Said, Steinert, and Mostajeran}

\ShortHeadings{Invariant kernels on the space of complex covariance matrices}{Said, Steinert, and Mostajeran}
\firstpageno{1}

\begin{document}

\title{Invariant Kernels on the space of \\ Complex Covariance Matrices}

\author{\name Salem Said \email salem.said@univ-grenoble-alpes.fr \\
       \addr Laboratoire Jean Kuntzman\\
       Universit\'e Grenoble-Alpes\\
       Grenoble, 38400, France \\
       \AND
       \name Franziskus Steinert \email franziskus.steinert@gmail.com \\
       \name Cyrus Mostajeran \email cyrussam.mostajeran@ntu.edu.sg \\
       \addr Division of Mathematical Sciences\\
		School of Physical and Mathematical Sciences\\
       Nanyang Technological University\\
       21 Nanyang Link, 637371, Singapore
       }

\editor{}

\maketitle

\begin{abstract}
  The present work develops certain analytical tools required to construct and compute invariant kernels on the space of complex covariance matrices.
  The main result is the $\mathrm{L}^1$--\,Godement theorem, which states that any invariant kernel, which is
  (in a certain natural sense) also integrable, can be computed by taking the inverse spherical transform of a positive function. General expressions for inverse spherical transforms are then provided, which can be used to explore new families of invariant kernels, at a rather moderate computational cost. A further, alternative approach for constructing new invariant kernels is also introduced, based on Ramanujan's master theorem for symmetric cones. This~leads to a novel closed-form invariant kernel, called the Beta-prime kernel. Numerical experiments highlight the computational and performance advantages of this kernel, especially in the context of two-sample hypothesis testing.  
\end{abstract}

\begin{keywords}
    kernel method, covariance matrix, Bochner's theorem, symmetric space,  
  spherical transform, Ramanujan's master theorem, two-sample test
\end{keywords}

\tableofcontents

\section{Introduction}
Positive definite kernels, and especially invariant positive definite kernels, play a prominent role across probability, statistics, and machine learning. Bochner's theorem characterizes the class of invariant positive definite kernels or (equivalently) the class of positive definite functions\,: a function $f:\mathbb{R}^{\scriptscriptstyle N} \rightarrow \mathbb{C}$ is positive definite if and only if it~is~the inverse Fourier transform of some finite positive measure. 

Recent research in machine learning and data science has focused on invariant positive definite kernels which~are defined on non-Euclidean spaces (rather than on $\mathbb{R}^{\scriptscriptstyle N}$ as in Bochner's theorem)~\citep{ruskernelspaper,ruskernelscompact,kernelspaper,rkks}. In turn, this revived interest in non-Euclidean generalizations of Bochner's theorem, due to Godement, Gelfand, and others~\citep{godement,yaglom}. In~\citep{kernelspaper}, the $\mathrm{L}^p$--\,Godement theorems (where $p = 1,2$) were introduced, in order to deal with integrable invariant kernels defined on non-compact Riemannian symmetric spaces. These are somewhat similar to Bochner's theorem, with the role of the Fourier transform played by the so-called spherical transform. 

The aim of the present work is to apply the $\mathrm{L}^1$--\,Godement theorem, in order to construct and compute positive definite functions on the space of complex covariance matrices. Our contribution is to develop the analytical results and tools that make it possible to work with spherical transforms on this space. We derive general expressions for inverse spherical transforms, which can be used to systematically explore new families of invariant positive definite functions with remarkable efficiency. We also develop an alternative approach for obtaining invariant positive definite functions in the form of spherical power series, based on Ramanujan's master theorem for symmetric cones. This is used to construct a novel closed-form invariant positive definite kernel with highly favorable computational properties, which we highlight through numerical examples. The results represent a significant advance in the theory and computation of positive definite kernels on spaces of covariance matrices.

Necessary background material is covered in Section \ref{sec:background}. This begins with two definitions, invariant positive definite kernels and $U$-invariant positive definite functions. Proposition \ref{prop:kerneltofunction} shows that these two concepts are in fact equivalent (below,~for brevity, the term "kernel" will mean positive definite kernel). Section \ref{sec:background} also discusses the Riemannian geometry of the space of complex covariance matrices, which is indeed a Riemannian symmetric space.

Section \ref{sec:sphericalt} is concerned with the spherical transform on this space. Roughly, this is an expansion of integrable $U$-invariant functions in terms of eigenfunctions of the Laplace-Beltrami operator, called spherical functions~\citep{helgasonbis,faraut}.\hfill\linebreak
Proposition \ref{prop:gn} provides a closed-form determinantal formula for these spherical functions, a slight generalization of the classic Gelfand-Naimark formula~\citep{gelfandnaimark}, while Proposition \ref{prop:spherical_pair} expresses the spherical transform itself, along with its inverse, in the form of a multiple integral. The following Proposition \ref{prop:gauss_spherical_trans} is the first application of Proposition~\ref{prop:spherical_pair}. It computes explicitly the spherical transform of a Gaussian function. 

Section \ref{sec:godement} comes to the main issue of applying the $\mathrm{L}^1$--\,Godement theorem. In~\citep{kernelspaper}, this theorem was stated under 
restrictive assumptions, which are not satisfied in the present work. Theorem \ref{th:L1G}, along with its proof in Appendix \ref{app:L1G}, make up for this problem. Proposition \ref{prop:langtrick} shows how Theorem \ref{th:L1G} can be used to construct and compute $U$-invariant positive definite functions. Specifically, the integrable $U$-invariant positive definite functions are exactly the inverse spherical transforms of a certain class of positive functions. Formulas (\ref{eq:PD1}) and (\ref{eq:PD2}) provide general expressions for these inverse spherical transforms. To generate a new positive definite function, it is enough to plug a suitable positive function into one of these two formulas. 

Section \ref{sec:godement} closes with three examples. The first one showcases an application of 
Formula (\ref{eq:PD1}). The second one uses (\ref{eq:PD2}) in order to compute the heat kernel (on the space of complex covariance matrices) in closed form. The last one combines Proposition \ref{prop:gauss_spherical_trans} and Theorem \ref{th:L1G} to prove that a Gaussian function cannot be positive definite.  

Section \ref{sec:ramanujan} features Proposition \ref{prop:ramanujan}, which shows how $U$-invariant positive definite functions can be obtained in the form of spherical power series, through an application of Ramanujan's master theorem for symmetric cones~\citep{ding}. An example of a positive definite function obtained in this way is the so-called "Beta-prime" function (see Expression (\ref{eq:betaprime})). The corresponding kernel will be called the Beta-prime kernel.

Section \ref{sec:application} is concerned with numerical experiments and applications. It begins with a computational comparison of the Beta-prime kernel to 
a prominent class of kernels from the literature, the heat and Matérn kernels \citep{ruskernelspaper}. It~then turns to the application of kernel methods to two-sample hypothesis testing, highlighting the computational advantages and performance capabilities of the Beta-prime kernel. 

It is remarkable that several spaces of covariance matrices (real, quaternion, block-Toeplitz, \textit{etc.}) can be embedded into spaces of complex covariance matrices (eventually of larger matrix size), in a way which preserves all the fundamental symmetry and invariance properties. Therefore, any invariant kernel constructed using the methods proposed in the present work immediately yields an invariant kernel on any of the above-mentioned spaces. In~this~way, the contribution of the present work is to introduce a general means of generating invariant kernels which admit analytical expressions and can be used on most of the usual spaces of covariance matrices. 

\section{General background} \label{sec:background} 

\subsection{Positive definite functions} \label{subsec:pdk}
Denote by $M$ the space of $N \times N$ complex covariance matrices. Specifically, these are $N \times N$ Hermitian positive definite matrices. Moreover, denote by $G$ the group of $N \times N$ invertible complex matrices, and by $U$ the group of $N \times N$ unitary matrices, a subgroup of $G$. 

Recall that  $G$ acts transitively on $M$ in the following way\,: $g\cdot x = gxg^\dagger$ for $g \in G$ and $x \in M$ (where $\dagger$ denotes the conjugate-transpose)~\citep{guigui,hdr}. For this action, $U$ is the stabiliser of the identity matrix $\mathrm{id} \in M$. In other words, $g \cdot \mathrm{id} = \mathrm{id}$ if and only if $g \in U$. 

A kernel $\mathcal{K}$ is a continuous function $\mathcal{K}: M \times M \rightarrow \mathbb{C}$, such that for any $x_{\scriptscriptstyle 1},\ldots,x_n \in M$ (here, $n = 2,3,\ldots$), the $n \times n$ matrix with elements $\mathcal{K}(x_i,x_j)$ is Hermitian positive semidefinite. The focus of the present work is on invariant kernels. These are kernels which satisfy $\mathcal{K}(g\cdot x,g \cdot y) = \mathcal{K}(x,y)$ for all $g \in G$ and all $x,y \in M$~\citep{ruskernelspaper,kernelspaper}. 

It is convenient to study invariant kernels indirectly, by studying $U$-invariant positive definite functions~\citep{kernelspaper}. A~function $f:M \rightarrow \mathbb{C}$ is called $U$-invariant if $f(u\cdot x) = f(x)$ for all $u \in U$ and $x \in M$. This means that $f(x) = f_o(\rho)$ where $f_o$ is a symmetric function and $\rho = (\rho_{\hspace{0.02cm}\scriptscriptstyle 1},\ldots,\rho_{\scriptscriptstyle N})$ are the eigenvalues of $x$\,: $f_o(\rho)$ remains unchanged after any permutation of $(\rho_{\hspace{0.02cm}\scriptscriptstyle 1},\ldots,\rho_{\scriptscriptstyle N})$.

If $f$ is continuous, then it is called positive definite if, for any $x_{\scriptscriptstyle 1},\ldots,x_n \in M$ (where $n = 2,3,\ldots$), the $n \times n$ matrix with elements $f(x^{\scriptscriptstyle -1/2}_ix_jx^{\scriptscriptstyle -1/2}_i)$ is Hermitian positive semidefinite. The two concepts (invariant kernel and $U$-invariant positive definite function) are equivalent. 
\begin{proposition} \label{prop:kerneltofunction}
If $\mathcal{K}$ is an invariant kernel, then the function $f(x) = \mathcal{K}(x,\mathrm{id})$ is $U$-invariant and positive definite. Conversely, if $f$ is a $U$-invariant and positive definite function, then $\mathcal{K}(x,y) = f(y^{\scriptscriptstyle -1/2}xy^{\scriptscriptstyle -1/2})$ defines an invariant kernel. 
\end{proposition}
\noindent \textbf{Proof\,:} let $f(x) = \mathcal{K}(x,\mathrm{id})$. If $u \in U$, then $f(u\cdot x) = \mathcal{K}(u\cdot x,\mathrm{id}) = \mathcal{K}(u\cdot x,u\cdot\mathrm{id})$ because $u\cdot \mathrm{id} = \mathrm{id}$.  However, if~$\mathcal{K}$ is invariant, then $\mathcal{K}(u\cdot x,u\cdot\mathrm{id})
= \mathcal{K}(x,\mathrm{id})$. Therefore, $f(u\cdot x) = f(x)$ and $f$ is $U$-invariant. To see that $f$ is~positive definite, it is enough to note that 
$$
f(y^{\scriptscriptstyle -1/2}xy^{\scriptscriptstyle -1/2}) = \mathcal{K}(y^{\scriptscriptstyle -1/2}\cdot x,\mathrm{id}) = 
\mathcal{K}(x,y^{\scriptscriptstyle 1/2}\cdot\mathrm{id}) = \mathcal{K}(x,y)
$$
where the first equality follows from the definition of $g\cdot x$, by taking $g = y^{\scriptscriptstyle -1/2}$, and the second equality because $\mathcal{K}$ is invariant. Thus, for any $x_{\scriptscriptstyle 1},\ldots,x_n \in M$, the matrix with elements $f(x^{\scriptscriptstyle -1/2}_ix_jx^{\scriptscriptstyle -1/2}_i)$ is~the same as the matrix with elements $\mathcal{K}(x_j,x_i)$, which is positive semidefinite because $\mathcal{K}$ is a kernel. In addition, continuity of $f$ follows from continuity of $\mathcal{K}$, and this ensures $f$ is positive definite.
 
Conversely, let $f$ be $U$-invariant and positive definite. This clearly implies that $\mathcal{K}(x,y) = f(y^{\scriptscriptstyle -1/2}xy^{\scriptscriptstyle -1/2})$ 
is~a~kernel. To see that this $\mathcal{K}$ is invariant, note that $\mathcal{K}(x,y) = f_o(\rho)$ where $\rho = (\rho_{\hspace{0.02cm}\scriptscriptstyle 1},\ldots,\rho_{\scriptscriptstyle N})$ are the eigenvalues of $y^{\scriptscriptstyle -1/2}xy^{\scriptscriptstyle -1/2}$. These are the same as the eigenvalues of $y^{\scriptscriptstyle -1}x$ because the two matrices are similar\,: $y^{\scriptscriptstyle -1}x = y^{\scriptscriptstyle -1/2}(y^{\scriptscriptstyle -1/2}xy^{\scriptscriptstyle -1/2})\hspace{0.03cm}y^{\scriptscriptstyle 1/2}\hspace{0.03cm}$. By the same argument, $\mathcal{K}(g\cdot x,g\cdot y) = f_o(\rho^\prime)$ where $\rho^\prime = (\rho^\prime_{\scriptscriptstyle 1},\ldots,\rho^\prime_{\scriptscriptstyle N})$ are the eigenvalues of $(g\cdot y)^{\scriptscriptstyle -1}(g\cdot x)$.  However, this last matrix is similar to $y^{\scriptscriptstyle -1}x$,
$$
(g\cdot y)^{\scriptscriptstyle -1}(g\cdot x) = (g^\dagger)^{\scriptscriptstyle -1}(y^{\scriptscriptstyle -1}x)(g^\dagger)
$$
Therefore, $\rho^\prime = \rho$ and $\mathcal{K}(g\cdot x,g\cdot y) = \mathcal{K}(x,y)$, as required. \hfill$\blacksquare$

\subsection{Riemannian geometry} \label{subsec:riemann}
An explicit description of $U$-invariant positive definite functions relies on the Riemannian geometry of $M$. Precisely, it relies on the fact that $M$ is a Riemannian symmetric space \citep{faraut,helgason}. 

Note that $M$ is an open subset of $H$, the real vector space of $N \times N$ Hermitian matrices. Therefore, $M$ is a differentiable manifold with its tangent space at any $x \in M$ naturally isomorphic to $H$. Now, with this in mind, consider the Riemannian metric on $M$, 
\begin{equation} \label{eq:aimetric}
  \langle v,w\rangle_x = \mathrm{Re}\left[\mathrm{tr}(x^{\scriptscriptstyle -1}v\hspace{0.02cm}x^{\scriptscriptstyle -1}w)\right] \hspace{1cm} v,w \in H
\end{equation}
where $\mathrm{Re}$ denotes the real part and $\mathrm{tr}$ the trace. This is the affine-invariant metric, which was made popular by~\citep{pennec}. For any $g \in G$, the map $x \mapsto g\cdot x$ is an isometry of the metric (\ref{eq:aimetric}). The same is true for inversion $x \mapsto x^{-1}$. These two facts together show that $M$ satisfies the definition of a Riemannian symmetric space. 

The class of $U$-invariant functions behaves in a special way with respect to the metric (\ref{eq:aimetric}). For example, let $\mathrm{vol}$ denote the Riemannian volume element of this metric. If $f:M \rightarrow \mathbb{C}$ is an integrable $U$-invariant function~\citep{hdr},
\begin{equation} \label{eq:uinv_integral}
   \int_M\,f(x)\hspace{0.02cm}\mathrm{vol}(dx) = \frac{C_{\scriptscriptstyle N}}{N!}\int_{\mathbb{R}^{\scriptscriptstyle N}_+}\,f_o(\rho)\hspace{0.02cm}(V(\rho))^2\hspace{0.02cm} \prod^N_{k=1} \rho^{-\scriptscriptstyle N}_kd\rho_k
\end{equation}
where $f(x) = f_o(\rho)$ is a symmetric function of the eigenvalues $(\rho_{\hspace{0.02cm}\scriptscriptstyle 1},\ldots,\rho_{\scriptscriptstyle N})$ of $x$, and where $V$ stands for the~Vandermonde polynomial. Here, and throughout the following, $C_{\scriptscriptstyle N}$ denotes a positive constant that only depends on $N$ and whose value is allowed to differ from one formula to another.

Moreover, if $L$ is the Laplace-Beltrami operator of the metric (\ref{eq:aimetric}), and $f$ is smooth and $U$-invariant, then~\citep{helgasonbis,faraut},  
\begin{equation} \label{eq:uinv_laplace}
   Lf = \sum^N_{k=1}\rho^2_k\frac{\partial^2f_o}{\mathstrut\partial \rho^2_k} + 2 \sum_{k<\ell}\frac{\rho_k\rho_\ell}{\rho_k - \rho_{\ell}}\left( \frac{\partial f_o}{\partial \rho_k} - \frac{\partial f_o}{\partial \rho_\ell}\right) + N \sum^N_{k=1}\rho_k\frac{\partial f_o}{\partial \rho_k}
\end{equation}
Formulas (\ref{eq:uinv_integral}) and (\ref{eq:uinv_laplace}) arise systematically from the Riemannian geometry of $M$, but they are also familiar in certain problems of multivariate statistics and random matrix theory~\citep{muirhead,mergny,mezzardi}.
       
To close the present section, consider a special case of the integral formula (\ref{eq:uinv_integral}). Assume that $f_o(\rho)$ factors into $f_o(\rho) = w(\rho_{\hspace{0.02cm}\scriptscriptstyle 1})\ldots w(\rho_{\scriptscriptstyle N})$ where $w$ is an integrable function such that 
$$
\int^\infty_0\,|w(\rho)|\hspace{0.03cm}\rho^{- N}\hspace{0.03cm}d\rho\,<\infty \hspace{0.25cm}\;\,\text{and}\hspace{0.25cm}
\int^\infty_0\,|w(\rho)|\hspace{0.03cm}\rho^{\hspace{0.03cm}N-2}\hspace{0.03cm}d\rho\,<\infty
$$
Then, the volume integral in (\ref{eq:uinv_integral}) is convergent and admits a determinantal expression
\begin{equation} \label{eq:determinantal_1}
   \int_M\,f(x)\hspace{0.02cm}\mathrm{vol}(dx) = C_{\scriptscriptstyle N}\,\det\!\left[\int^\infty_0\,w(\rho)\hspace{0.02cm}\rho^{k+\ell- N}d\rho\right]^{N-1}_{k,\ell=0}   
\end{equation}
This is an application of the Andr\'eief identity, widely used in random matrix theory, and was pointed out in~\citep{tierz}. \\[0.1cm]
\textbf{Example\,:} the expression (\ref{eq:determinantal_1}) can be used to compute the Gaussian integral of~\citep{said2,said3}
\begin{equation} \label{eq:gauss_Z}
  Z(\sigma) = \int_M\,\exp\left[-\frac{d^{\hspace{0.03cm} 2}(x,\mathrm{id})}{2\sigma^2}\right]\mathrm{vol}(dx)
\end{equation}
where $d(\cdot,\cdot)$ denotes the Riemannian distance induced on $M$ by the metric (\ref{eq:aimetric}). In (\ref{eq:determinantal_1}), this integral corresponds to $w(\rho) = \exp[-\log^2(\rho)/2\sigma^2]$, which yields
\begin{equation} \label{eq:gauss_det}
  Z(\sigma) = C_{\scriptscriptstyle N}\,\det\!\left[\sigma\hspace{0.02cm}e^{(\sigma^2/2)(k+\ell-N-1)^2}\right]^{N}_{k,\ell=1}
\end{equation}
a formula due to~\citep{tierz}, which will be both simplified and generalized in the following section (see Proposition \ref{prop:gauss_spherical_trans}). 

\section{The spherical transform} \label{sec:sphericalt}
The key ingredient which will be employed in constructing and computing $U$-invariant positive definite functions is the spherical transform. Roughly, this provides an expansion of any well-behaved $U$-invariant function, in terms of eigenfunctions of the Laplace-Beltrami operator (\ref{eq:uinv_laplace}), which are known as spherical functions~\citep{helgasonbis,faraut}. 

The set of all spherical functions is described as follows~\citep{faraut}. Consider first the power function, $\Delta_s : M \rightarrow \mathbb{C}$ where $s = (s_{\scriptscriptstyle 1},\ldots, s_{\scriptscriptstyle N})$ belongs to $\mathbb{C}^{\scriptscriptstyle N}$, 
\begin{equation} \label{eq:powerfunction}
\Delta_s(x) = (\Delta_{\scriptscriptstyle 1}(x))^{ s_{\scriptscriptstyle 1}-s_{\scriptscriptstyle 2}}(\Delta_{\scriptscriptstyle 2}(x))^{ s_{\scriptscriptstyle 2} - s_{\scriptscriptstyle 3}}\ldots (\Delta_{\scriptscriptstyle N}(x))^{ s_{\scriptscriptstyle N}}
\end{equation}
where $\Delta_k(x)$ is the $k$-th leading principal minor of $x \in M$. A spherical function is a function of the form 
\begin{equation} \label{eq:sphericalfunction}
 \Phi_\lambda(x) = \int_U\, \Delta_{\lambda + \delta}(u\cdot x)\hspace{0.03cm}du
\end{equation}
where $\lambda \in \mathbb{C}^{\scriptscriptstyle N}$ and $\delta_k = \frac{1}{2}(2k-N-1)$, while $du$ denotes the normalized Haar measure on the unitary group $U$. Two functions $\Phi_{\lambda}$ and $\Phi_{\lambda^\prime}$ are identical if and only if $(\lambda_{\scriptscriptstyle 1},\ldots, \lambda_{\scriptscriptstyle N})$ is a permutation of $(\lambda^\prime_{\scriptscriptstyle 1},\ldots, \lambda^\prime_{\scriptscriptstyle N})$. 

Each $\Phi_\lambda$ is $U$-invariant and an eigenfunction of the Laplace-Beltrami operator (\ref{eq:uinv_laplace}), with eigenvalue $(\lambda,\lambda) - (\delta,\delta)$ \citep{faraut}~(Theorem XIV.3.1). Here, and throughout the following, $(\mu\hspace{0.02cm},\nu) = \sum^N_{k=1} \mu_k\hspace{0.02cm}\nu_k$ for $\mu,\nu\in \mathbb{C}^{\scriptscriptstyle N}$. 

The first claim of the present section is that the spherical functions $\Phi_\lambda\hspace{0.03cm}$, while initially given by the integral formula (\ref{eq:sphericalfunction}), admit the following determinantal expression. 
\begin{proposition} \label{prop:gn}
   If  $\lambda \in \mathbb{C}^{\scriptscriptstyle N}$ and $\Phi_\lambda$ is given by (\ref{eq:sphericalfunction}), then
\begin{equation} \label{eq:gelfandnaimark}
 \Phi_\lambda(x) = \prod^{N-1}_{k=1}k!\times
\frac{\det\!\left[\rho^{\lambda_{\ell}+\scriptscriptstyle (N-1)/2}_k\right]^{N}_{k,\ell=1}}{V(\lambda)V(\rho)}
\end{equation} 
where $(\rho_{\hspace{0.02cm}\scriptscriptstyle 1},\ldots, \rho_{\scriptscriptstyle N})$ are the eigenvalues of $x$ and $V$ stands for the~Vandermonde polynomial. 
\end{proposition}
 
The proof of Proposition \ref{prop:gn} will be given in Appendix \ref{app:gn}. Formula (\ref{eq:gelfandnaimark}) will be called the Gelfand-Naimark formula, as it is a slight generalization of the formula introduced by Gelfand and Naimark~\citep{gelfandnaimark}. The reason why the spherical functions (\ref{eq:sphericalfunction}) admit the determinantal expression (\ref{eq:gelfandnaimark}) is that the group $G$ is here a complex Lie group. This~fact is the foundation of the proof in Appendix \ref{app:gn}.

Now, let $f:M \rightarrow \mathbb{C}$ be an integrable $U$-invariant function (integrable means with respect to $\mathrm{vol}$, as in (\ref{eq:uinv_integral})).
Its~spherical transform is the function $\hat{f}:\mathbb{R}^{\scriptscriptstyle N} \rightarrow \mathbb{C}$, 
\begin{equation} \label{eq:sphericaltransform}
  \hat{f}(t) = \int_M\,f(x)\Phi_{-\mathrm{i}t}(x)\hspace{0.02cm}\mathrm{vol}(dx)
\end{equation}
where $\mathrm{i} = \sqrt{-1}$~\citep{faraut}. An inversion theorem for the spherical transform (\ref{eq:sphericaltransform}) is given in~\citep{faraut} (Theorem XIV.5.3). Specifically, if $\hat{f}$ satisfies the integrability condition
\begin{equation} \label{eq:hc_integrable}
  \int_{\mathbb{R}^{\scriptscriptstyle N}}\,|\hat{f}(t)|\hspace{0.03cm}(V(t))^2dt < \infty
\end{equation}
then the following inversion formula holds,
\begin{equation} \label{eq:spherical_inverse}
   f(x) = C_{\scriptscriptstyle N}\int_{\mathbb{R}^{\scriptscriptstyle N}}\,\hat{f}(t)\Phi_{\hspace{0.03cm}\mathrm{i}t}(x)\hspace{0.02cm}(V(t))^2dt
\end{equation}
After substituting (\ref{eq:uinv_integral}) and (\ref{eq:gelfandnaimark}) into (\ref{eq:sphericaltransform}) and (\ref{eq:spherical_inverse}), the following is obtained.
\begin{proposition} \label{prop:spherical_pair}
  For the spherical transform pair (\ref{eq:sphericaltransform})--(\ref{eq:spherical_inverse}),
\begin{equation} \label{eq:sphericaltransform_bis}
  \hat{f}(t) = \frac{C_{\scriptscriptstyle N}}{V(-\mathrm{i}t)} \times \frac{1}{N!} \int_{\mathbb{R}^{\scriptscriptstyle N}_+}\,f_o(\rho)\hspace{0.02cm}V(\rho)\det\!\left[\rho^{-\mathrm{i}t_{\ell}-\scriptscriptstyle (N+1)/2}_k\right]\hspace{0.02cm}d\rho
\end{equation}
\begin{equation} \label{eq:spherical_inverse_bis}
  f(x) =  \frac{C_{\scriptscriptstyle N}}{V(\mathrm{i}\rho)}\times 
\frac{1}{N!}\int_{\mathbb{R}^{\scriptscriptstyle N}}\,\hat{f}(t)
V(t)\det\!\left[\rho^{\mathrm{i}t_{\ell}+\scriptscriptstyle (N-1)/2}_k\right]\hspace{0.02cm}dt
\end{equation}
where $f(x) = f_o(\rho)$ is a symmetric function of the eigenvalues $(\rho_{\scriptscriptstyle 1},\ldots,\rho_{\scriptscriptstyle N})$ of $x$.    
\end{proposition}

The proof of this proposition will not be given in detail, as it merely consists of performing straightforward algebraic simplifications.

It is remarkable that the spherical transform does not involve all the spherical functions $\Phi_\lambda$ but only the functions $\Phi_{\hspace{0.03cm}\mathrm{i}t}$ where $t \in \mathbb{R}^{\scriptscriptstyle N}$. These functions have in common the property that they correspond to real, negative eigenvalues of the Laplace-Beltrami operator (\ref{eq:uinv_laplace})\,: each $\Phi_{\hspace{0.03cm}\mathrm{i}t}$ corresponds to the eigenvalue $-(t,t) - (\delta,\delta)$. The spherical functions that do not appear in the spherical transform are interesting in their own right. For instance, if $\lambda + \delta = m$ where $(m_{\scriptscriptstyle 1},\ldots,m_{\scriptscriptstyle N})$ are positive integers arranged in decreasing order, then one has
\begin{equation} \label{eq:schur}
\Phi_\lambda(x) = \prod^{N-1}_{k=1}k!\times\left.S_m(\rho)\middle/V(\lambda)\right.
\end{equation}
where $S_m$ denotes the Schur polynomial corresponding to $(m_{\scriptscriptstyle 1},\ldots,m_{\scriptscriptstyle N})$. Schur polynomials are very important in the study of circular and unitary-invariant random matrix ensembles~\citep{mergny,mezzardi,meckes} (because they provide the irreducible characters of the unitary group). In Section \ref{sec:ramanujan}, below, they will appear within the framework of Ramanujan's master theorem. 

The following proposition is motivated by the study of the Gaussian integral (\ref{eq:gauss_Z}). Specifically, consider the integrals
\begin{equation} \label{eq:characteristic}
  Z(\sigma,\lambda) = \int_M\,\exp\left[-\frac{d^{\hspace{0.02cm}\scriptscriptstyle 2}(x,\mathrm{id})}{2\sigma^2}\right]\Phi_{\lambda}(x)\,\mathrm{vol}(dx)
\end{equation}
where $\Phi_\lambda$ was defined in (\ref{eq:sphericalfunction}). If $\lambda = -\delta$ then $Z(\sigma,\lambda)$ is just $Z(\sigma)$ from (\ref{eq:gauss_Z}). On the other hand, note that $Z(\sigma,-\mathrm{i}t) = \hat{f}(t)$ where $f(x) = \exp[-d^{\hspace{0.02cm}\scriptscriptstyle 2}(x,\mathrm{id})/2\sigma^2]$. 
\begin{proposition} \label{prop:gauss_spherical_trans}
    The integrals (\ref{eq:characteristic}) admit the following expression
 \begin{flalign} 
\nonumber & Z(\sigma,\lambda) = 
 \frac{C_{\scriptscriptstyle N}}{V(\lambda)} \times \det\!\left[\,\sigma\hspace{0.03cm}\exp\left((\sigma^2/2)\!\left(\delta_k+\lambda_\ell\right)^{2}\hspace{0.03cm}\right)\right]^N_{k,\ell=1} \\[0.15cm]
\label{eq:gauss_ZT} &  = C_{\scriptscriptstyle N}\hspace{0.03cm}\sigma^{\scriptscriptstyle N^2}\, e^{\frac{\sigma^2}{2}\left((\lambda,\lambda)+(\delta,\delta)\right)}\,\prod_{k < \ell}\mathrm{sch}\!\left((\sigma^2/2)(\lambda_\ell - \lambda_k)\right)
 \end{flalign} 
 where $\mathrm{sch}(a) = \sinh(a)/a$.
\end{proposition}
The proof of Proposition \ref{prop:gauss_spherical_trans} will be given in Appendix \ref{app:gauss_st}. This proposition provides a simplified form of (\ref{eq:gauss_det}),
\begin{equation} \label{eq:gauss_Z_new}
   Z(\sigma) = 
   C_{\scriptscriptstyle N}\hspace{0.03cm}\sigma^{\scriptscriptstyle N}\, e^{\sigma^2(\delta,\delta)}\,\prod_{k < \ell}\sinh\!\left((\sigma^2/2)(\ell - k)\right)
\end{equation}
as follows by replacing $\lambda = -\delta$ into (\ref{eq:gauss_ZT}). In addition, for the Gaussian function $f$ defined before the proposition, putting $\lambda = -\mathrm{i}t$ gives the spherical transform
\begin{equation} \label{eq:gauss_spherical_transform}
  \hat{f}(t) = C_{\scriptscriptstyle N}\hspace{0.03cm}\sigma^{\scriptscriptstyle N^2}\, e^{\frac{\sigma^2}{2}\left((\delta,\delta)-(t,t)\right)}\,\prod_{k < \ell}\mathrm{sc}\!\left((\sigma^2/2)(t_\ell - t_k)\right)
\end{equation}
where $\mathrm{sc}(a) = \sin(a)/a$. 

\section{Constructing invariant kernels} \label{sec:godement}
The $\mathrm{L}^1$--\,Godement theorem was introduced in~\citep{kernelspaper}. Roughly, this theorem shows that $U$-invariant positive definite functions can be obtained by taking inverse spherical transforms of positive symmetric functions. 
\begin{theorem} \label{th:L1G}
    Let $f:M \rightarrow \mathbb{C}$ be an integrable $U$-invariant function (integrable means with respect to $\mathrm{vol}$).
    Then, $f$ is positive definite if and only if
\begin{equation} \label{eq:L1G}
  f(x) = \int_{\mathbb{R}^{\scriptscriptstyle N}}\,g(t)\Phi_{\hspace{0.03cm}\mathrm{i}t}(x)\hspace{0.02cm}(V(t))^2dt
\end{equation}
where the function $g:\mathbb{R}^{\scriptscriptstyle N}\rightarrow \mathbb{R}$ is positive ($g(t) \geq 0$ for all $t$), symmetric, and satisfies the integrability condition
\begin{equation} \label{eq:hc_integrable_1}
  \int_{\mathbb{R}^{\scriptscriptstyle N}}\,g(t)\hspace{0.03cm}(V(t))^2dt < \infty
\end{equation}
Moreover, this function $g$ is then unique --- recall that $g$ is said to be symmetric if $g(t)$ remains unchanged after any permutation of $(t_{\scriptscriptstyle 1},\ldots,t_{\scriptscriptstyle N})$.
\end{theorem}

The proof of Theorem \ref{th:L1G} will be given in Appendix \ref{app:L1G}. As explained in~\citep{kernelspaper}, the $\mathrm{L}^1$--\,Godement theorem is based on the celebrated Godement theorem, which generalizes Bochner's theorem to the context of symmetric spaces~\citep{godement}. The only-if part of this theorem can be used to check whether a given function $f$ is positive definite or not. On the other hand, the if part can be used to construct and compute positive definite functions. Indeed, note that (\ref{eq:L1G}) is essentially an inverse spherical transform as in (\ref{eq:spherical_inverse}), with $g(t)$ instead of $\hat{f}(t)$. Therefore, just as in (\ref{eq:spherical_inverse_bis}), it is possible to rewrite (\ref{eq:L1G}),
\begin{equation} \label{eq:L1G_bis}
   f(x) =  \frac{1}{V(\mathrm{i}\rho)}\times 
\frac{1}{N!}\int_{\mathbb{R}^{\scriptscriptstyle N}}\,g(t)
V(t)\det\!\left[\rho^{\mathrm{i}t_{\ell}+\scriptscriptstyle (N-1)/2}_k\right]\hspace{0.02cm}dt
\end{equation}
where $(\rho_{\scriptscriptstyle 1},\ldots,\rho_{\scriptscriptstyle N})$ are the eigenvalues of $x$. To obtain a positive definite function $f$, it is then enough to choose a suitable positive function $g$ and then evaluate the integral (\ref{eq:L1G_bis}). 


\begin{proposition} \label{prop:langtrick}
  Let $g:\mathbb{R}^{\scriptscriptstyle N}\rightarrow \mathbb{R}$ satisfy the conditions of Theorem \ref{th:L1G}. 
\begin{enumerate} 
\item[(a)] Assume that $g(t)$ factors into $g(t) = \gamma(t_{\scriptscriptstyle 1})\ldots\gamma(t_{\scriptscriptstyle N})$, where $\gamma$ is a positive function. It~follows from (\ref{eq:L1G_bis}) that
\begin{equation} \label{eq:PD1}
    f(x) = \frac{(\det(x))^{\scriptscriptstyle {(N-1)/2}}}{V(\mathrm{i}\rho)}\times \det\!\left[\int_{\mathbb{R}}\,\gamma(t)t^{k-1}\hspace{0.03cm}e^{\mathrm{i}t s_\ell}\hspace{0.03cm}dt \right] \\[0.12cm]
\end{equation}
Here, $s_\ell = \log(\rho_\ell)$ for $\ell = 1,\ldots, N$. 
\item[(b)] Assume that the inverse Fourier transform
$$
\tilde{g}(s) = \int_{\mathbb{R}^{\scriptscriptstyle N}}\,g(t)\hspace{0.03cm}e^{\mathrm{i}(s,t)}\hspace{0.02cm}dt \hspace{1cm} (s_{\scriptscriptstyle 1},\ldots, s_{\scriptscriptstyle N}) \in \mathbb{R}^{\scriptscriptstyle N}
$$
is smooth. It follows from (\ref{eq:L1G_bis}) that
\begin{equation} \label{eq:PD2}
    f(x) = \frac{(\det(x))^{\scriptscriptstyle {(N-1)/2}}}{V(\rho)}\times \left.V\left(-\frac{\partial}{\partial s}\right)\tilde{g}(s)\right|_{s_\ell =
 \log(\rho_\ell)} \\[0.12cm]
\end{equation}
where $V(\partial/\partial s)$ is the Vandermonde operator $V(\partial/\partial s) = \prod_{i<j} (\partial/\partial s_j - \partial/\partial s_i)$. 
\end{enumerate}
\end{proposition}
The proof of Proposition \ref{prop:langtrick} is given in Appendix \ref{app:langtrick}. \\

\noindent \textbf{Example\,:} as an application of Proposition \ref{prop:langtrick}-(a), choose $g(t) = \gamma(t_{\scriptscriptstyle 1})\ldots\gamma(t_{\scriptscriptstyle N})$ where $\gamma(t) = (\kappa/2)\exp(-\kappa|t|)$ for~$\kappa > 0$.  Replacing into (\ref{eq:PD1}) and using elementary properties of the Fourier transform,
\begin{equation} \label{eq:cauchy}
  f(x) = \frac{(\det(x))^{\scriptscriptstyle {(N-1)/2}}}{V(\rho)}\times \det\!\left[-\tilde{\gamma}^{(k-1)}(\log(\rho_\ell))\right] 
\end{equation}
where $\tilde{\gamma}^{(k-1)}$ is the $(k-1)$-th derivative of $\tilde{\gamma}(s) = (\kappa^2+s^2)^{-1}\hspace{0.03cm}$. Theorem \ref{th:L1G} now says that this $f$ is a positive definite function. In fact, (\ref{eq:cauchy}) provides a whole family of positive definite functions, parameterized by $\kappa$. \\

\noindent \textbf{Example\,:} while it looks a bit too complicated, Proposition \ref{prop:langtrick}-(b) can be used to compute the heat kernel of $M$ in closed form~\citep{helgasonbis,lang}. This corresponds to $g(t) = \exp[-\kappa((t,t)+(\delta,\delta))]$ with $\kappa > 0$. Now, for this choice of $g(t)$, one has the following inverse Fourier transform
$$
\tilde{g}(s) = \left((\pi/\kappa)^{\scriptscriptstyle N/2}\hspace{0.03cm}e^{-\kappa(\delta,\delta)}\right)\exp\left[-(s,s)/4\kappa\right]
$$
and one may use a beautiful identity from~\citep{lang} (Chapter XII, Page 405)
\begin{equation} \label{eq:langidentity}
 V\left(-\frac{\partial}{\partial s}\right)\exp\left[-(s,s)/4\kappa\right] = (1/2\kappa)^{\scriptscriptstyle N(N-1)/2}\hspace{0.03cm}V(s)\exp\left[-(s,s)/4\kappa\right]
\end{equation}
Replacing this into (\ref{eq:PD2}) yields the heat kernel (rather, $f(x) = \mathcal{K}(x,\mathrm{id})$ where $\mathcal{K}$ is the heat kernel)
\begin{equation} \label{eq:heatkernel}
  f(x) = C_\kappa\hspace{0.03cm}(\det(x))^{\scriptscriptstyle {(N-1)/2}}\times \left(V(\log(\rho))\middle/V(\rho)\right)\exp\left[-|\log(\rho)|^2/4\kappa\right]
\end{equation}
where $C_\kappa > 0$ is a constant and $|\log(\rho)|^2 = (\log(\rho),\log(\rho))$. Of course, this $f$ is a positive definite function for~each $\kappa > 0$, thanks to Theorem \ref{th:L1G}. \\

\noindent \textbf{Example\,:} recall the Gaussian function $f$ defined before Proposition \ref{prop:gauss_spherical_trans}, with its spherical transform $\hat{f}$ in (\ref{eq:gauss_spherical_transform}). Unlike the functions in the two previous examples, this one is never positive definite. Precisely, there exists no value of $\sigma$ for which it is positive definite. This~can be seen using the only-if part of Theorem \ref{th:L1G}. According to this theorem, since $f$ is integrable and $U$-invariant, if $f$ were positive definite, its spherical transform $\hat{f}$ would be identical to the positive function $g$ in (\ref{eq:L1G}), up~to a constant factor. Indeed, because the spherical transform $\hat{f}$ in (\ref{eq:gauss_spherical_transform}) satisfies the integrability condition (\ref{eq:hc_integrable}), the inversion formula (\ref{eq:spherical_inverse}) holds true. The fact that $g = C_{\scriptscriptstyle N}\hspace{0.02cm}\hat{f}$ then follows by injectivity of the inverse spherical transform. 
Now, to show that $f$ is not positive definite, it is enough to show that $\hat{f}(t) < 0$ for some $t \in \mathbb{R}^{\scriptscriptstyle N}$. Choosing $t$ such that $t_{\scriptscriptstyle 1} < \ldots < t_{\scriptscriptstyle N}\hspace{0.03cm}$, note that $t_{\ell} - t_k \leq t_{\scriptscriptstyle N} - t_{\scriptscriptstyle 1}$ for any $k < \ell$, with equality only if $(k,\ell) = (1,N)$. If $t_{\scriptscriptstyle N} - t_{\scriptscriptstyle 1} \leq \pi$ then $\hat{f}(t) \geq 0$, but as soon as $t_{\scriptscriptstyle N} - t_{\scriptscriptstyle 1} > \pi$, then $\mathrm{sc}\!\left((\sigma^2/2)(t_{\scriptscriptstyle N} - t_{\scriptscriptstyle 1})\right)$ becomes the only negative term in the product on the right-hand side of (\ref{eq:gauss_spherical_transform}), and it then follows that $\hat{f}(t) < 0$. 

\section{Applying Ramanujan's theorem} \label{sec:ramanujan}
Ramanujan's master theorem for symmetric cones was stated and proved in~\citep{ding}. Roughly, this theorem converts so-called spherical power series into spherical transforms (or inverse spherical transforms). It is a generalization of the theorem obtained by Ramanujan, which deals with classical (one-variable) power series and Mellin transforms~\citep{hardy}. Here, the aim is to examine how this theorem can be used in order to construct $U$-invariant positive definite functions. In the first place, this is possible because the space $M$ of $N \times N$ complex covariance matrices is indeed a symmetric cone~\citep{faraut}. Spherical power series will be expressed in terms of complex zonal polynomials~\citep{takemura}. The Schur polynomials $S_m$ from (\ref{eq:schur}) in Section \ref{sec:sphericalt} can be normalized to~obtain new polynomials $Z_m$ with the property that 
$$
(\mathrm{tr}(x))^n = \sum_{[m]=n} Z_m(\rho) \hspace{0.5cm} \text{where } 
[m] = m_{\scriptscriptstyle 1}+\ldots+m_{\scriptscriptstyle N}
$$
for any positive integer $n$. These are the complex zonal polynomials, and a series of the following form
\begin{equation} \label{eq:spherical_taylor}
    F(x) = \sum_m \frac{(-1)^{[m]}}{[m]!}\hspace{0.02cm}a(m)\hspace{0.03cm}Z_m(x)
\end{equation}
will be called a spherical power series (the sum is over positive integers $m_{\scriptscriptstyle 1} \geq \ldots \geq m_{\scriptscriptstyle N} \geq 0$). 

Ramanujan's master theorem can be used to prove the following Proposition \ref{prop:ramanujan}. In this proposition, $\Gamma_M$ denotes the Gamma function of the symmetric cone $M$,
\begin{equation} \label{eq:GammaM}
\Gamma_M(\lambda) = (2\pi)^{N(N-1)/2}\prod^N_{k=1} \Gamma(\lambda_k - k + 1) \hspace{1cm} \lambda \in \mathbb{C}^{\scriptscriptstyle N}
\end{equation}
where $\Gamma$ is Euler's Gamma function of a complex variable~\citep{faraut,ding}. Moreover, the variable $x$ in (\ref{eq:spherical_taylor}) ranges over the real vector space $H$ of $N \times N$ Hermitian matrices ($M$ is an open cone within $H$). 
\begin{proposition} \label{prop:ramanujan}
Let $\alpha > N-1$ and assume that the coefficients $a(m)$ in (\ref{eq:spherical_taylor}) are given by a function $a:\mathbb{C}^{\scriptscriptstyle N} \rightarrow \mathbb{C}$, of the form
\begin{equation} \label{eq:rama_cond1}
a(\lambda) = \Gamma_M(2\alpha+\lambda)\hspace{0.02cm}\psi(\lambda - \delta)
\end{equation}
where $\psi:\mathbb{C}^{\scriptscriptstyle N} \rightarrow \mathbb{C}$ is symmetric, holomorphic for $\mathrm{Re}(\lambda_k) > (N-1)/2 - 2\alpha$, and satisfies the growth condition
\begin{equation} \label{eq:rama_cond2}
    \left| \psi(\lambda)\right| \leq C_{\scriptscriptstyle N}\hspace{0.02cm}\prod^N_{k=1} e^{P\hspace{0.02cm}\mathrm{Re}(\lambda_k)}\times 
    e^{A\hspace{0.02cm}|\mathrm{Im}(\lambda_k)|}
\end{equation}
where $P,A > 0$ and $A <\pi$. Then, the series (\ref{eq:spherical_taylor}) converges in a neighborhood of $x = 0$, where it defines a real-analytic function $F$. Moreover, this function extends analytically to all of $M$, in such a way that $\Delta^\alpha(x)F(x)$ defines an integrable $U$-invariant function of $x \in M$ ($\Delta(x) = \det(x))$, which is positive definite if and only if $\psi(\mathrm{i}t - \alpha) \geq 0$ for all $t \in \mathbb{R}^{\scriptscriptstyle N}$. 
\end{proposition}

The proof of Proposition \ref{prop:ramanujan} is given in Appendix \ref{app:ramanujan}. The idea is to show that $f(x) = \Delta^\alpha(x)F(x)$ defines an integrable $U$-invariant function with spherical transform
\begin{equation} \label{eq:ramanujan}
  \hat{f}(t) = \left|\Gamma_M(\alpha + \delta + \mathrm{i}t)\right|^2\psi(\mathrm{i}t - \alpha)
\end{equation}
and that this $\hat{f}(t)$ satisfies the integrability condition (\ref{eq:hc_integrable}), so that the inversion formula (\ref{eq:spherical_inverse}) is also satisfied. The~statement about positive definiteness of $f$ then follows by an application of Theorem \ref{th:L1G}. 

\begin{figure}
\centering
\includegraphics[width=0.75\linewidth]{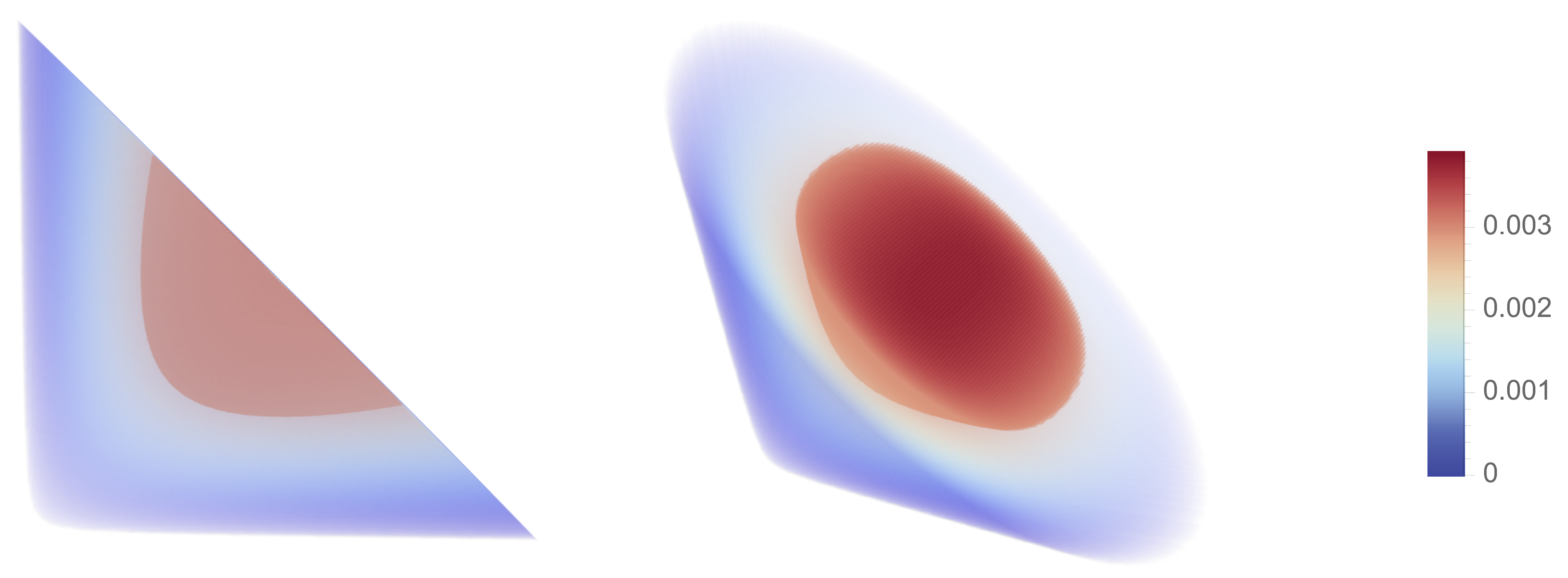}
\caption{Color plot of $f(x)/\Gamma_M(2\alpha)$, given by (\ref{eq:betaprime}) with $\alpha = 2$, restricted to real $2\times 2$ SPD matrices, here visualized as points in a cone viewed from two perspectives. The domain of the depicted plot is restricted to $0<\operatorname{tr}(x)<2$.}\label{fig:betaprime_2x2}
\end{figure}

\noindent \textbf{Example:} in (\ref{eq:rama_cond1}), if $\psi$ is a constant function equal to $1$, the series (\ref{eq:spherical_taylor}) becomes a generalized binomial series~\citep{faraut} (Proposition XII.1.3), which converges to $F(x) = \Gamma_M(2\alpha)\Delta^{-2\alpha}(\mathrm{id}+x)$ for any $x$ whose eigenvalues are all $< 1$. Proposition \ref{prop:ramanujan} then says that, for each $\alpha > N-1$,
\begin{equation} \label{eq:betaprime}
 f(x) = \Gamma_M(2\alpha)\hspace{0.02cm}\Delta^\alpha\left(x\hspace{0.02cm}(\mathrm{id}+x)^{-2}\right)
\end{equation}
is an integrable $U$-invariant positive definite function, whose spherical transform can be read from (\ref{eq:ramanujan}), by setting $\psi \equiv 1$. In the one-dimensional case ($N = 1$), this function reduces to the density of the Beta-prime distribution~\citep{kotz}. For $N = 2$, Figure \ref{fig:betaprime_2x2} shows a color plot of $f(x)/\Gamma_M(2\alpha)$ with $\alpha = 2$. In this plot, $x$ is restricted to real, symmetric positive definite (SPD) matrices. A $2\times 2$ SPD matrix $x=\left[x_{ij}\right]_{i,j=1,2}$ is represented as a point $(x_{11},x_{12},x_{22}) \in \mathbb{R}^3$, which lies in the cone described by $x_{11}>0$ and $x_{11}x_{22}-x_{12}^2 > 0$.

\section{Numerical Experiments} \label{sec:application}

\subsection{Computational considerations}
According to Proposition \ref{prop:kerneltofunction}, the function $f$ in (\ref{eq:betaprime}) corresponds to an invariant kernel $\mathcal{K}(x,y) = f(y^{\scriptscriptstyle -1/2}xy^{\scriptscriptstyle -1/2})$. After some elementary simplifications, this can be written
\begin{equation} \label{eq:betaprime_k}
    \mathcal{K}(x,y) = \Gamma_M(2\alpha)\, \bigg[ \frac{\Delta(x)\, \Delta(y)}{\Delta(x+y)^2}\bigg]^\alpha\hspace{1cm} \alpha > N-1
\end{equation}
where $\Delta(x) = \det(x)$, just as in (\ref{eq:betaprime}). This kernel will be called the \emph{Beta-prime kernel}. 

Evaluation of this kernel does not involve matrix multiplication or inversion, only computation of~determinants, which is typically implemented via LU-factorization. For $N \times N$ matrices, it has computational complexity $\mathcal{O}(N^3)$. The factor $\Gamma_M(2\alpha)$ ahead of (\ref{eq:betaprime_k}) may of course be dropped without affecting positive-definiteness of the kernel. 

In this section, the Beta-prime kernel will be considered as a kernel on the space of real covariance matrices. This~only involves requiring $x$ and $y$ in (\ref{eq:betaprime_k}) to be real, and therefore symmetric positive-definite (SPD) matrices. In~\citep{ruskernelspaper}, alternative kernels on the space of real SPD matrices have been considered\,: the heat and Matérn kernels.
The striking difference between the Beta-prime kernel and these alternative kernels is the significantly lower computational complexity of the Beta-prime kernel.

In fact, the heat and Matérn kernels on the space of real SPD matrices do not have closed-form expressions. In~\citep{ruskernelspaper,Mostowsky2024geometrickernelspackage}, they are  approximated using a dedicated Monte Carlo technique. On the space of $N \times N$ SPD matrices, this is used to evaluate a rather complicated $N$-dimensional integral, by means of a double Monte Carlo integration, in addition to a rejection sampling step. It should then come as no surprise the closed-form Beta-prime kernel has a low complexity in comparison to such an elaborate approach.
A further drawback of Monte Carlo techniques is that they introduce noise in the kernel values, which may even lead to loss of positive definiteness.

\begin{figure}
\centering
\includegraphics[width=0.49\linewidth]{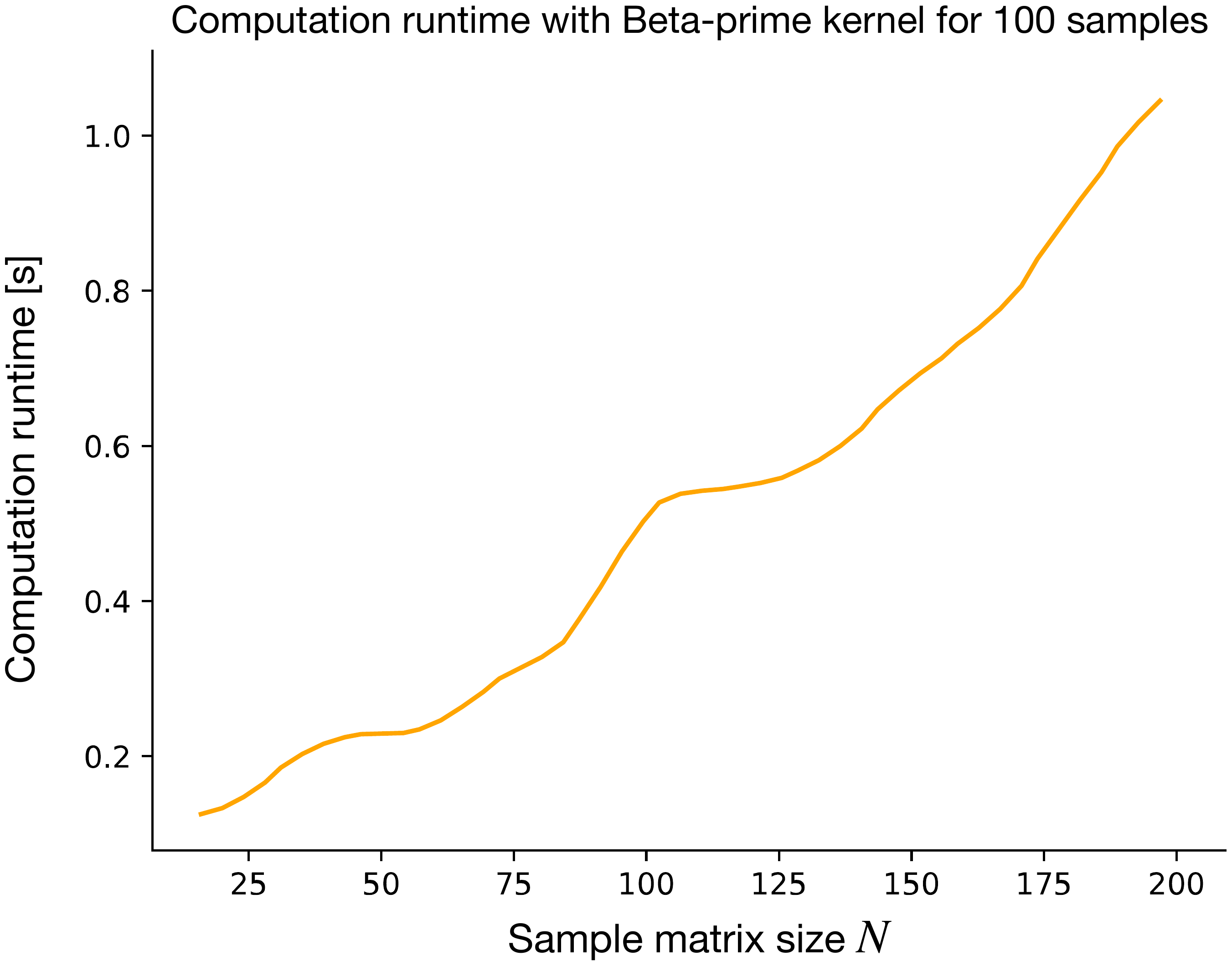}
\includegraphics[width=0.49\linewidth]{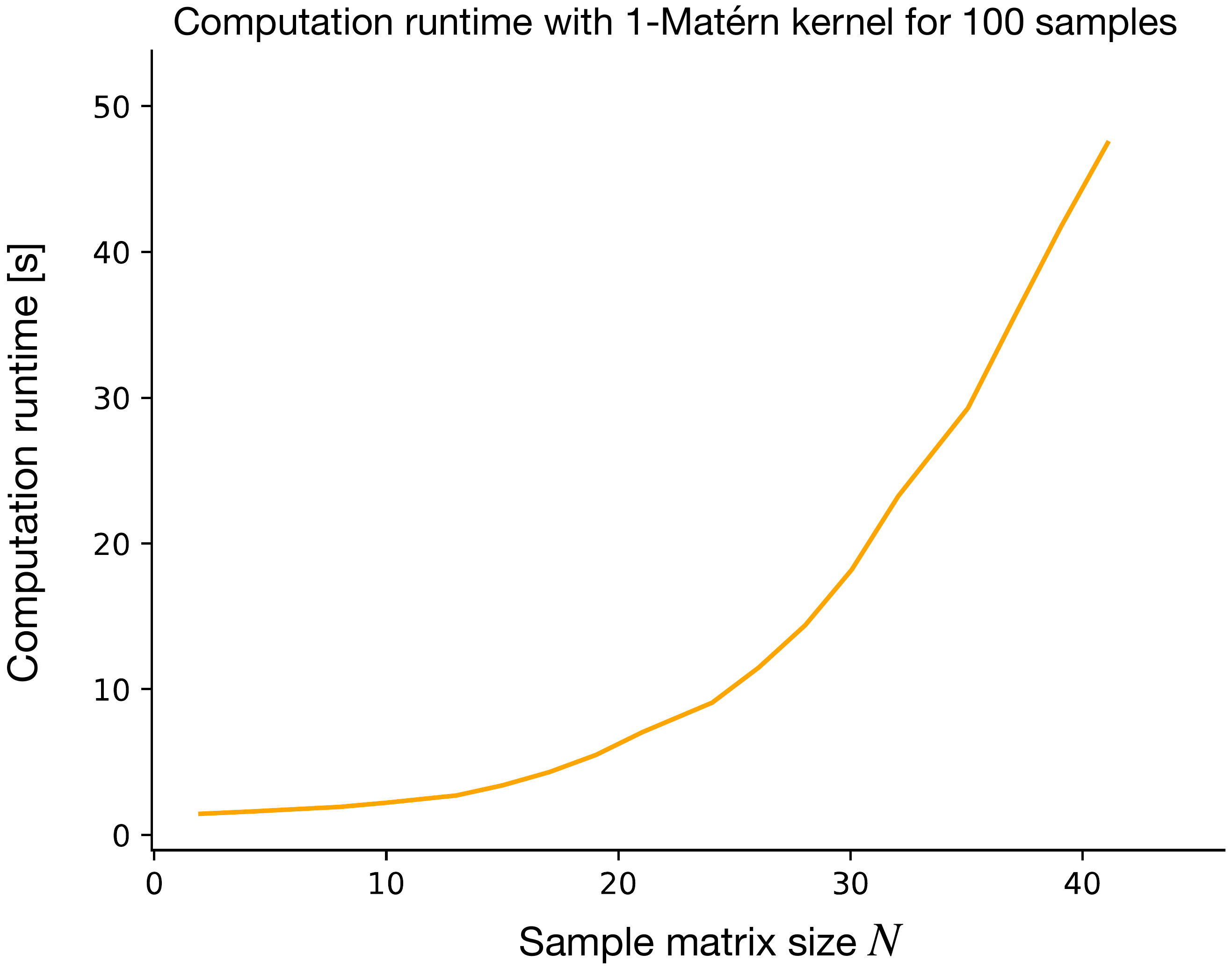}
\caption{Computation runtime for $100 \times 100$ kernel Gram matrices, plotted against matrix dimension $N$. Beta-prime kernel (left) and 1-Matérn kernel (right). The $1$-Matérn kernel is the Matérn kernel with smoothness parameter $ = 1$.} \label{fig:scaling_dim}
\end{figure}

To visualize the difference in computational complexity, Figure~\ref{fig:scaling_dim} shows
runtime (in seconds) as a function of the matrix dimension $N$. Specifically, the runtime is the average time required to compute the kernel Gram matrix $[\mathcal{K}(x_i,x_j)]^{100}_{i,j=1}$ of a random sample of $N \times N$ SPD matrices $(x_i;i=1,\ldots, 100)$. 

Especially notable is the scale difference on the axes of the two plots. The plot on the left involves the Beta-prime kernel from (\ref{eq:betaprime_k}) with $\alpha = N$. The plot on the right involves the Matérn kernel with smoothness parameter equal to $1$~\citep{ruskernelspaper} (see Proposition 12).

\subsection{Application: the two-sample test}
A fundamental problem in statistics is to determine whether two given samples $X = \{ x_1, \ldots , x_m\} \sim p$ and 
$Y=\{ y_1, \ldots ,y_{m'}\} \sim q$ originate from the same distribution on a space $\mathcal{X}$. This is the two-sample test problem. Two-sample tests are important to many applications, 
and are used to compare the results of generative models~\citep{lopez2016revisiting}.

Several hypothesis testing methods have been developed that summarize the difference in the statistical properties of $X$ and $Y$ in a test statistic, which is then used to decide whether to reject the null hypothesis "$p=q$" \citep{lopez2016revisiting}.
A~prominent test statistic in this context is maximum mean discrepancy (MMD),
\begin{equation} \label{MMD}
    MMD_\mathcal{F}\left(p,q\right):= \sup_{f\in \mathcal{F}} |\mathbb{E}_{x\sim p}(f(x))-\mathbb{E}_{y\sim q}(f(y))|
\end{equation}
where $\mathcal{F}$ is a certain class of functions $f:\mathcal{X}\rightarrow\mathbb{R}$.
This quantity, of course, depends on the choice of $\mathcal{F}$. When~choosing $\mathcal{F}=\{ f\in \mathcal{C}_b(\mathcal{X}): |f|\leq 1\}$, where $\mathcal{C}_b(\mathcal{X})$ is the space of bounded continuous functions on $\mathcal{X}$, one~has $p=q$ if and only if $MMD_\mathcal{F}\left(p,q\right)=0$. However, one may also choose a smaller $\mathcal{F}$ that is dense in $\mathcal{C}_b(\mathcal{X})$ with respect to the uniform topology. Such $\mathcal{F}$ is provided by the unit ball in the RKHS of a $\mathcal{C}$-universal kernel $\mathcal{K}$.
In this case~\citep{gretton2012kernel}, the squared MMD statistic is computed as 
\begin{equation} \label{MMDkernel}
    MMD^2_\mathcal{F}\left(p,q\right)= \mathbb{E}_{x,x'}[\mathcal{K}(x,x')]-2\, \mathbb{E}_{x,y}[\mathcal{K}(x,y)] + \mathbb{E}_{y,y'}[\mathcal{K}(y,y')]
\end{equation}
where $x'$ and $y'$ are independent copies of $x$ and $y$.

Assuming that the samples $X$ and $Y$ are of equal size $m$, an unbiased estimate of $MMD^{2}_\mathcal{F}$ was given by~\citep{gretton2012kernel},
\begin{equation} \label{MMDunbiased}
    \widehat{MMD}^2_\mathcal{F}\left(p,q\right) = \frac{1}{m\,(m-1)} \sum_{i\neq j}h(z_i,z_j)
\end{equation}
where $z_i=(x_i, y_i)$ and 
\begin{equation*}
    h(z_i, z_j) = \mathcal{K}(x_i, x_j)- \mathcal{K}(x_i,y_j)-\mathcal{K}(x_j,y_i)+ \mathcal{K}(y_i, y_j)
\end{equation*}
For the present application, $\mathcal{X} = \mathcal{P}_N$ is the space of $N \times N$ SPD matrices. The Beta-prime kernel from (\ref{eq:betaprime_k}) can be proved to be universal for continuous functions on compact subsets of $\mathcal{P}_N$ \citep{Steinert2025}. Hence, if the samples $X,Y\subset \mathcal{P}_N$ are drawn from distributions supported in a compact subset of $\mathcal{P}_N\hspace{0.03cm}$, one may use the Beta-prime kernel to implement the empirical test statistic (\ref{MMDunbiased}). The same is true of the heat and Matérn kernels on $\mathcal{P}_N\hspace{0.03cm}$, but these have to be computed using the computationally costly approach of~\citep{ruskernelspaper}.

Several two-sample tests were performed, with both kinds of kernels (Beta-prime and Matérn). The matrix dimension was $N = 3$ and the sample size $m = 100$. For $X$, eigenvectors were sampled from a uniform distribution on the three-dimensional orthogonal group, while eigenvalues were independent, uniformly distributed in the interval $[30,31]$. Several new samples $Y_k$ had the same eigenvector distribution, but eigenvalues uniformly distributed in
$[r_k\cdot 30, \, r_k\cdot 30+1]$, where the \emph{spectral scaling factors} $r_k$ were given by $r_k= 0.1+k(4-0.1)/80$ for $k=0,\ldots, 80$.

While the tests using the empirical test statistic (\ref{MMDunbiased}) were completed within seconds using the Beta-prime kernel, those using the $1$-Matérn kernel were prohibitively slow. 
All tests were therefore redone, using the linear time test statistic from~\citep{gretton2012kernel} (Lemma 14). The aim was to determine whether $X$ and $Y_k$ were drawn from the same distribution on a significance level of 5\%.

\begin{figure}
\centering
\includegraphics[width=0.9\linewidth]{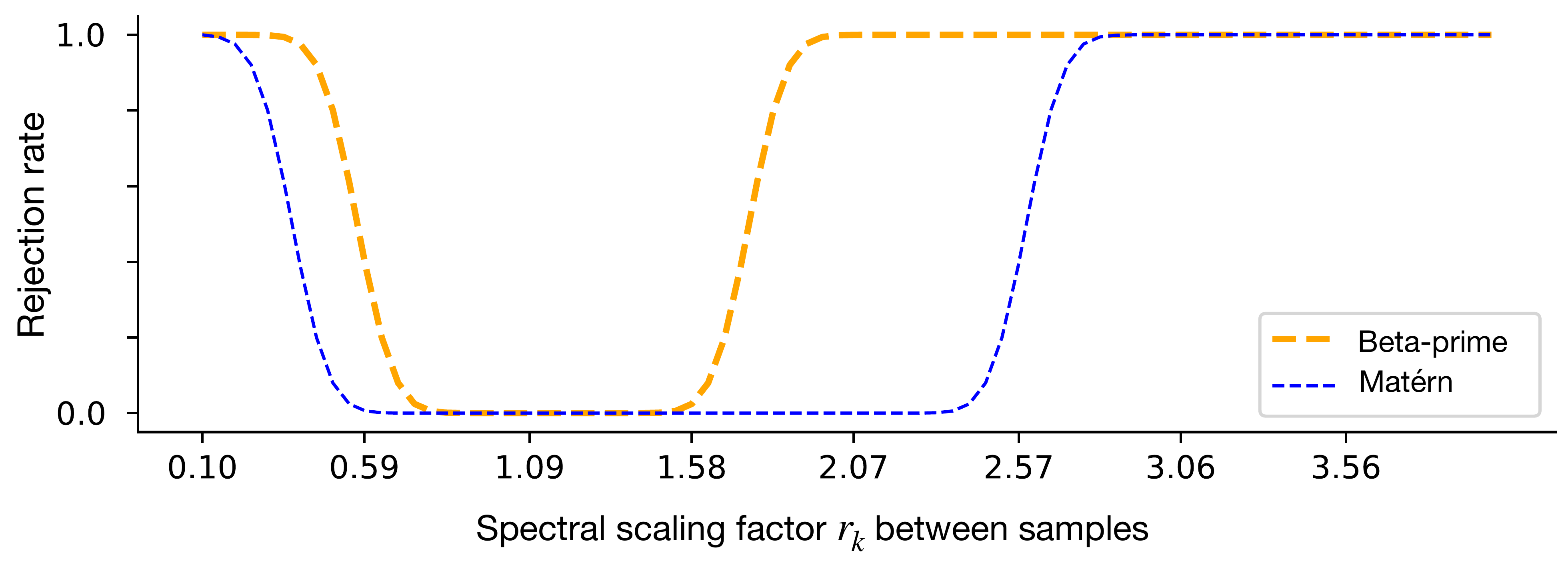}
\caption{Rejection rate of the null hypothesis "same distribution" plotted against the spectral scaling factor $r_k\hspace{0.03cm}$.
Beta-prime kernel (yellow) and 1-Matérn kernel (blue).}
\label{fig:two_sample_test}
\end{figure}

The results are illustrated in Figure~\ref{fig:two_sample_test} which shows the relative number of times that the null hypothesis "same distribution" was rejected (rejection rate) for different values of $r_k\hspace{0.03cm}$. Tests using the Beta-prime kernel were completed in approximately one second on a personal laptop computer, whereas the same tests using the $1$-Matérn kernel required approximately eleven hours (using the approach of~\citep{ruskernelspaper} and the bespoke computational package~\citep{ Mostowsky2024geometrickernelspackage}). 

In~addition to this dramatic computational advantage, the Beta-prime kernel exhibited a clear advantage in performance over the $1$-Matérn kernel, as can be seen in Figure~\ref{fig:two_sample_test}, where it is clearly much more effective at distinguishing the two distributions. Indeed, the range of values of $r_k$ where the rejection rate is low (the dip in the curves in Figure \ref{fig:two_sample_test}) is narrower and more concentrated near $r_k = 1$ for the yellow (Beta-prime) curve than for the blue (Matérn) curve. 

\acks{The authors acknowledge financial support from the School of Physical and Mathematical Sciences and the Talent Recruitment and Career Support (TRACS) Office at Nanyang Technological University (NTU Singapore).}
    
\vskip 0.2in
\bibliography{references}

\vfill
\pagebreak

\appendix

\section{Proof of Proposition \ref{prop:gn}} \label{app:gn}
The proof starts from (\ref{eq:sphericalfunction}) and shows that it is equivalent to (\ref{eq:gelfandnaimark}). The first step is to reduce to the case where $x$ has unit determinant. Let $x = r\hspace{0.02cm}\bar{x}$ where $\det(\bar{x}) = 1$. By a direct calculation, it follows from
(\ref{eq:powerfunction}) that $\Delta_s(x) = r^{(1,s)}\Delta_s(\bar{x})$ where $(1,s) = \sum^N_{k=1} s_k\hspace{0.03cm}$. In particular, $s = \lambda + \delta$ implies that $\Delta_s(x) = r^{(1,\lambda)}\Delta_s(\bar{x})$, because $(1,\delta) = 0$. Replacing into (\ref{eq:sphericalfunction}) yields the following identity
\begin{equation} \label{eq:proof_gn1}
   \Phi_\lambda(x) = r^{(1,\lambda)}\hspace{0.02cm}\Phi_\lambda(\bar{x})
\end{equation}
Now, returning to (\ref{eq:sphericalfunction}), note that
\begin{equation} \label{eq:proof_gn2}
\Phi_\lambda(\bar{x}) = \int_U\, \Delta_{\lambda + \delta}(u\cdot \bar{x})\hspace{0.03cm}du
\end{equation}
This integral can be restricted to $SU$, the special unitary subgroup of $U$ ($SU$ is the set of $u \in U$ with $\det(u) = 1$). 
Indeed, replacing $u$ by $e^{\mathrm{i}\theta}\hspace{0.02cm}u$ (with $\theta$ real) does not change $u\cdot \bar{x}$. Moreover, the normalized Haar measure on $U$ descends to the normalized Haar measure on $SU$~\citep{knapp} (see Theorem 8.32). Therefore, (\ref{eq:proof_gn2}) is equivalent to
\begin{equation} \label{eq:proof_gn3}
\Phi_\lambda(\bar{x}) = \int_{SU}\, \Delta_{\lambda + \delta}(u\cdot \bar{x})\hspace{0.03cm}du
\end{equation}
The next step of the proof is to show that (\ref{eq:proof_gn3}) is the same as the following Harish-Chandra integral~\citep{helgasonbis} (Page 418)
\begin{equation} \label{eq:proof_gn4}
\Phi_\lambda(\bar{x}) = \int_{SU}\,\exp\left[\left(2\lambda-2\delta,\log a\!\left(u\bar{\rho}^{\scriptscriptstyle\hspace{0.03cm} 1/2}\hspace{0.02cm}\right)\right)\right]du
\end{equation}
Here, $\bar{\rho} = (\bar{\rho}_{\scriptscriptstyle 1},\ldots,\bar{\rho}_{\scriptscriptstyle N})$ are the eigenvalues of $\bar{x}$, and $a(u\bar{\rho}^{\hspace{0.02cm}\scriptscriptstyle 1/2})$ is the diagonal matrix with positive entries, such that $u\bar{\rho}^{\hspace{0.02cm}\scriptscriptstyle 1/2} = n\hspace{0.02cm}a(u\bar{\rho}^{\hspace{0.02cm}\scriptscriptstyle 1/2})\hspace{0.02cm}h$ where the matrix $n$ is upper triangular  with ones on its diagonal, and $h$ belongs to $SU$ (vectors such as $\bar{\rho}$ will be identified with diagonal matrices, in a self-evident way, whenever that is convenient). 

Because $\Phi_\lambda$ is $U$-invariant, it is enough to prove (\ref{eq:proof_gn4}) when $\bar{x} = \bar{\rho}$. It will be helpfull to apply the identity $\Phi_\lambda(x) = \Phi_{-\lambda}(x^{-\scriptscriptstyle 1})$~\citep{faraut} (Theorem XIV.3.1). For short, let $a = a(u\bar{\rho}^{\hspace{0.02cm}\scriptscriptstyle 1/2})$, so $u\cdot \bar{\rho} = (u\bar{\rho}^{\hspace{0.02cm}\scriptscriptstyle 1/2})\cdot \mathrm{id}$ is equal to 
$n\hspace{0.02cm}a^2\hspace{0.02cm}n^\dagger$. Using the identity just mentioned, 
\begin{equation} \label{eq:proof_gn5}
\Phi_\lambda(\bar{\rho}) = \Phi_{-\lambda}(\bar{\rho}^{\hspace{0.02cm}-\scriptscriptstyle 1}) = \int_{SU}\, \Delta_{-\lambda + \delta}\left(u\cdot \bar{\rho}^{\hspace{0.02cm}- 1}\right)\hspace{0.03cm}du 
=
\int_{SU}\,\Delta_{-\lambda + \delta}\left((n^\dagger)^{-1}\cdot a^{-2}\right)\hspace{0.03cm}du 
\end{equation}
where the second equality follows from (\ref{eq:proof_gn3}) and the third equality holds because $u\cdot \bar{\rho}^{\hspace{0.02cm}-\scriptscriptstyle 1} = (u\cdot \bar{\rho})^{-\scriptscriptstyle 1}$ and $u\cdot \bar{\rho} = n\hspace{0.02cm}a^2\hspace{0.02cm}n^\dagger$. 
However, since $(n^\dagger)^{-1}$ is lower triangular with ones on its diagonal, it is easy to see that $
\Delta_{-\lambda + \delta}\left((n^\dagger)^{-1}\cdot a^{-2}\right) = \Delta_{-\lambda + \delta}\left(a^{-2}\right)
$.
 Then, from (\ref{eq:powerfunction}) and the fact that $a$ is diagonal,
$$
\Delta_{-\lambda + \delta}\left((n^\dagger)^{-1}\cdot a^{-2}\right) = \prod^N_{k=1}\,a^{2(\lambda_k-\delta_k)}_k = 
\exp\left[\left(2\lambda-2\delta,\log a\right)\right]
$$
and (\ref{eq:proof_gn4}) follows immediately by replacing this into (\ref{eq:proof_gn5}).  The final step of the proof exploits the fact that $G$ is a complex Lie group. In this case, the Harish-Chandra integral (\ref{eq:proof_gn4}) admits a closed-form expression~\citep{helgasonbis} (Page 432),
\begin{equation} \label{eq:complexcase}
  \Phi_\lambda(\bar{x}) = \frac{\Pi(-\delta)}{\Pi(\lambda)}\times\frac{\sum_{w \in S_{\scriptscriptstyle N}} \varepsilon(w)e^{(\lambda,w \bar{\rho})}}{\sum_{w \in S_{\scriptscriptstyle N}}\varepsilon(w)e^{(\delta,w \bar{\rho})}}
\end{equation}
Here, $\Pi(\mu) = \prod_{k<\ell} (\mu_k - \mu_\ell)$ for $\mu\in \mathbb{C}^{\scriptscriptstyle N}$, $S_{\scriptscriptstyle N}$ is the symmetric group (group of permutations of $N$ objects), and $\varepsilon(w)$ is the signature of the permutation $w$, while $w \bar{\rho}$ denotes the action of that permutation on $(\bar{\rho}_{\scriptscriptstyle 1},\ldots,\bar{\rho}_{\scriptscriptstyle N})$. Clearly, the polynomial $\Pi$ is the Vandermonde polynomial, up to sign, while the sums in the second fraction are Leibniz expansions of certain determinants. Using these observations, and performing some basic simplifications,
\begin{equation} \label{eq:proof_gn6}
  \Phi_\lambda(\bar{x}) = \frac{\prod^{N-1}_{k=1}k!}{V(-\lambda)}\times\frac{\det\!\left[\bar{\rho}^{\hspace{0.02cm}\lambda_{\ell}}_k\right]}{\det\!\left[
\bar{\rho}^{\hspace{0.02cm}N-\ell}_k\hspace{0.03cm}\right]} =
\frac{\prod^{N-1}_{k=1}k!}{V(-\lambda)}\times\frac{\det\!\left[\bar{\rho}^{\hspace{0.02cm}\lambda_{\ell}}_k\right]}{V(-\bar{\rho})} 
= \prod^{N-1}_{k=1}k! \times \frac{\det\!\left[\bar{\rho}^{\hspace{0.02cm}\lambda_{\ell}}_k\right]}{V(\lambda)V(\bar{\rho})} 
\end{equation}
Now, (\ref{eq:gelfandnaimark}) can be retrieved from (\ref{eq:proof_gn1}) and (\ref{eq:proof_gn6}). To do so, it is enough to note that $\rho = r\hspace{0.02cm}\bar{\rho}$ and that this implies (recall $\rho = (\rho_{\scriptscriptstyle 1},\ldots,\rho_{\scriptscriptstyle N})$ are the eigenvalues of $x$)
$$
r^{(1,\lambda)}\hspace{0.02cm}\det\!\left[\bar{\rho}^{\hspace{0.02cm}\lambda_{\ell}}_k\right]  = \det\!\left[\rho^{\lambda_{\ell}}_k\right] 
\hspace{0.5cm}\text{and}\hspace{0.5cm} V(\bar{\rho}) = r^{-\scriptscriptstyle N(N-1)/2}\hspace{0.03cm}V(\rho) 
$$
Then, since $r^{\scriptscriptstyle N} = \det(x)$ (which is the product of the eigenvalues $\rho_k$),
\begin{equation} \label{eq:proof_gn7}
r^{(1,\lambda)}\hspace{0.02cm}\frac{\det\!\left[\bar{\rho}^{\hspace{0.02cm}\lambda_{\ell}}_k\right]}{V(\bar{\rho})} = r^{\scriptscriptstyle N(N-1)/2}\hspace{0.02cm}\frac{\det\!\left[\rho^{\lambda_{\ell}}_k\right]}{V(\rho)} = \frac{\det\!\left[\rho^{\lambda_{\ell}+\scriptscriptstyle (N-1)/2}_k\right]}{V(\rho)} 
\end{equation}
\vspace{0.2cm}
Therefore, performing the multiplication in (\ref{eq:proof_gn1}), with the help of (\ref{eq:proof_gn6}) and (\ref{eq:proof_gn7}), yields the required (\ref{eq:gelfandnaimark}).

\section{Proof of Proposition \ref{prop:gauss_spherical_trans}}\label{app:gauss_st}
Let $f(x) = \exp[-d^{\hspace{0.02cm}\scriptscriptstyle 2}(x,\mathrm{id})/2\sigma^2]$. Then, note that (\ref{eq:characteristic}) reads
\begin{equation} \label{eq:proofchar1}
 Z(\sigma,\lambda) = \int_M\,f(x)\hspace{0.03cm}\Phi_{\lambda}(x)\,\mathrm{vol}(dx)
\end{equation}
Since both $f$ and $\Phi_\lambda$ are $U$-invariant, this can be evaluated using (\ref{eq:uinv_integral}). In terms of the eigenvalues $(\rho_{\hspace{0.02cm}\scriptscriptstyle 1},\ldots,\rho_{\scriptscriptstyle N})$, $f(x) = w(\rho_{\hspace{0.02cm}\scriptscriptstyle 1})\ldots w(\rho_{\scriptscriptstyle N})$ where $w(\rho) = \exp[-\log^2(\rho)/2\sigma^2]$~\citep{tierz,said3}. On the other hand, $\Phi_\lambda(x)$ is given by (\ref{eq:gelfandnaimark}). Replacing this into (\ref{eq:uinv_integral}), it follows that
$$
Z(\sigma,\lambda) = 
\frac{C_{\scriptscriptstyle N}}{V(\lambda)} \times \frac{1}{N!} \int_{\mathbb{R}^{\scriptscriptstyle N}_+}\,V(\rho)\det\!\left[\rho^{\lambda_{\ell}-\scriptscriptstyle (N+1)/2}_k\right]\hspace{0.02cm}\prod_k w(\rho_k)d\rho_k
$$
Recalling that $V(\rho) = \det\left[\rho^{\ell-1}_k\right]$, and applying the Andr\'eief identity (as stated in~\citep{livan}, Chapter 11), this becomes
$$
Z(\sigma,\lambda) = 
\frac{C_{\scriptscriptstyle N}}{V(\lambda)} \times
\det\!\left[\int^\infty_0\,w(\rho)\hspace{0.02cm}
\rho^{k - 1 + \lambda_\ell - \scriptscriptstyle (N+1)/2} \hspace{0.03cm}d\rho\right]^N_{k,\ell=1}
$$
By the definition of $\delta_k$ (right after Formula (\ref{eq:sphericalfunction}) in Section \ref{sec:sphericalt}), this is the same as
$$
Z(\sigma,\lambda) = 
\frac{C_{\scriptscriptstyle N}}{V(\lambda)} \times
\det\!\left[\int^\infty_0\,w(\rho)\hspace{0.02cm}
\rho^{\delta_k + \lambda_\ell - 1} \hspace{0.03cm}d\rho\right]^N_{k,\ell=1}
$$
The integrals inside the determinant can be expressed in terms of the moments of a log-normal probability density. This gives the first line in (\ref{eq:gauss_ZT}),
\begin{equation} \label{eq:proof_char2}
Z(\sigma,\lambda) = 
 \frac{C_{\scriptscriptstyle N}}{V(\lambda)} \times \det\!\left[\,\sigma\hspace{0.03cm}\exp\left((\sigma^2/2)\!\left(\delta_k+\lambda_\ell\right)^{2}\hspace{0.03cm}\right)\right]^N_{k,\ell=1}
\end{equation}
To complete the proof of (\ref{eq:gauss_ZT}), it is enough to use elementary properties of the determinant,
\begin{flalign}
    \nonumber \det\!\left[\,\exp\left((\sigma^2/2)\!\left(\delta_k+\lambda_\ell\right)^{2}\hspace{0.03cm}\right)\right] & = \det\!\left[\,\exp\left((\sigma^2/2)\!\left(\delta^{\hspace{0.02cm} 2}_k+\lambda^2_\ell+2\delta_k\hspace{0.02cm}\lambda_\ell\right)\hspace{0.03cm}\right)\right] \\[0.15cm]
\label{eq:proof_char3}    &  = \prod^N_{k=1} e^{\frac{\sigma^2}{2}\left(\delta^{\hspace{0.02cm} 2}_k+\lambda^2_k\right)} \times 
\det\!\left[\,\exp
\left(\sigma^2\hspace{0.02cm}\delta_k\hspace{0.02cm}\lambda_\ell\right)\hspace{0.03cm}\right]
\end{flalign}
Indeed, it is clear that 
\begin{equation} \label{eq:proof_char4}
\prod^N_{k=1} e^{\frac{\sigma^2}{2}\left(\delta^{\hspace{0.02cm} 2}_k+\lambda^2_k\right)} = e^{\frac{\sigma^2}{2}\left((\lambda,\lambda)+(\delta,\delta)\right)}
\end{equation}
Moreover, using the definition of $\delta_k$ and performing some straightforward simplifications
\begin{flalign}
\nonumber \det\!\left[\,\exp
\left(\sigma^2\hspace{0.02cm}\delta_k\hspace{0.02cm}\lambda_\ell\right)\hspace{0.03cm}\right] & =
\prod^N_{k=1} e^{\frac{-\sigma^2}{2}(N-1)\lambda_k} \times V(e^{\sigma^2\hspace{0.02cm}\lambda}) \\[0.15cm]    
\label{eq:proof_char5} & = \prod_{k < \ell} 2\hspace{0.02cm}\sinh\!\left((\sigma^2/2)(\lambda_\ell - \lambda_k)\right)
\end{flalign}
Therefore, replacing (\ref{eq:proof_char4}) and (\ref{eq:proof_char5}) into (\ref{eq:proof_char3}), 
$$
\det\!\left[\,\exp\left((\sigma^2/2)\!\left(\delta_k+\lambda_\ell\right)^{2}\hspace{0.03cm}\right)\right] = e^{\frac{\sigma^2}{2}\left((\lambda,\lambda)+(\delta,\delta)\right)}\,\prod_{k < \ell} 2\hspace{0.02cm}\sinh\!\left((\sigma^2/2)(\lambda_\ell - \lambda_k)\right)
$$
In turn, replacing this into (\ref{eq:proof_char2}) directly yields the second line in (\ref{eq:gauss_ZT}).

\section{Proof of Theorem \ref{th:L1G}} \label{app:L1G}
\noindent \textbf{if part\,:} roughly, the idea of the proof is to show that the spherical functions $\Phi_{\hspace{0.03cm}\mathrm{i}t}$ are positive definite functions. Then, (\ref{eq:L1G}) says that the function $f$ is a positive linear combination of these $\Phi_{\hspace{0.03cm}\mathrm{i}t}$ and is therefore positive definite. For $x \in M$, let $x = \exp(s)\hspace{0.03cm}\bar{x}$ where $s = \log\!\hspace{0.02cm}\det(x)$ so that $\det(\bar{x}) = 1$. From (\ref{eq:proof_gn1}) (putting $\lambda = \mathrm{i}\hspace{0.02cm}t$ and $\tau = (1,t)$),
\begin{equation} \label{eq:proof_L1G_if_1}
  \Phi_{\hspace{0.03cm}\mathrm{i}t}(x) = e^{\mathrm{i}\tau s}\hspace{0.03cm} \Phi_{\hspace{0.03cm}\mathrm{i}t}(\bar{x})
\end{equation}
Thinking of $s$ and $\bar{x}$ as functions of $x$, let $\varphi^\tau(x) = e^{\mathrm{i}\tau s}$ and $\bar{\varphi}(x) = \Phi_{\hspace{0.03cm}\mathrm{i}t}(\bar{x})$, so that $\Phi_{\hspace{0.03cm}\mathrm{i}t}(x) = \varphi^\tau(x)\hspace{0.03cm}\bar{\varphi}(x)$. Now, recalling the well-known fact that a product of positive definite functions is positive definite, it is enough to show that $\varphi^\tau$ and $\bar{\varphi}$ are both positive definite. For any $x_{\scriptscriptstyle 1},\ldots,x_n \in M$, note that
\begin{equation} \label{eq:PDmatrix1}
\varphi^\tau(x^{\scriptscriptstyle -1/2}_ix_jx^{\scriptscriptstyle -1/2}_i) = e^{\mathrm{i}\tau (s_j - s_i)}
\end{equation}
where $s_i = \log\!\hspace{0.02cm}\det(x_i)$. Therefore, the matrix with elements $\varphi^\tau(x^{\scriptscriptstyle -1/2}_ix_jx^{\scriptscriptstyle -1/2}_i)$ is the same as the matrix with elements $e^{\mathrm{i}\tau (s_j - s_i)}$, which is Hermitian non-negative definite (of rank $1$). This shows that $\varphi^\tau$ is positive definite. To see that $\bar{\varphi}$ is also positive definite, note that
\begin{equation} \label{eq:PDmatrix2}
\bar{\varphi}(x^{\scriptscriptstyle -1/2}_ix_jx^{\scriptscriptstyle -1/2}_i) = \Phi_{\hspace{0.03cm}\mathrm{i}t}(\bar{x}^{\scriptscriptstyle -1/2}_i\bar{x}_j\bar{x}^{\scriptscriptstyle -1/2}_i)
\end{equation}
However, according to~\citep{helgasonbis} (Page 484), the restriction of $\Phi_{\hspace{0.03cm}\mathrm{i}t}$ to the unit-determinant hypersurface (that is to the set of $x \in M$ with $\det(x) = 1$), which is given by the Harish-Chandra integral (\ref{eq:proof_gn4}), is a positive definite function. In particular, the matrix whose elements appear in (\ref{eq:PDmatrix2}) is Hermitian non-negative definite. This shows that $\bar{\varphi}$ is positive definite. Thus, being a product of positive definite functions, $\Phi_{\hspace{0.03cm}\mathrm{i}t}$ is positive definite (for any $t$) as required.\hfill\linebreak
To conclude, recall Godement's theorem~\citep{kernelspaper,godement} (in particular, Formula (2.4) in~\citep{kernelspaper}). This gives a rigorous justification of the claim that $f$ is positive definite because it is a positive linear combination of positive definite functions. 
The~last step of the proof is thus a direct application of Godement's theorem. 

\noindent \textbf{only-if part\,:} as stated in~\citep{kernelspaper}, the $\mathrm{L}^1$--\,Godement theorem says that an integrable function $f$ is positive definite if~and only if its spherical transform is positive and integrable (details can be found in~\citep{kernelspaper}, Section 2 and Appendix A).\hfill\linebreak 
In~\citep{kernelspaper}, it is required from the outset that the underlying symmetric space $M$ should be of non-compact type (in~particular, the group $G$ should be semisimple). This requirement is clearly not satisfied, in the present context.  
The aim here is to explain that the $\mathrm{L}^1$--\,Godement theorem can still be applied.

In~\citep{kernelspaper}, the requirement that $M$ should be of non-compact type was introduced only in order to ensure that spherical functions on $M$ are given by Harish-Chandra integrals (integrals of the form (\ref{eq:proof_gn4}) in the proof of Proposition \ref{prop:gn}). 
The proof of the $\mathrm{L}^1$--\,Godement theorem (see~\citep{kernelspaper}, Appendix A) relies on Godement's (much older) theorem~\citep{godement}, which applies to any symmetric space and in particular to the space $M$ of complex covariance matrices. Specifically, the convolution product of two compactly-supported continuous $U$-invariant functions on $M$ is commutative~\citep{faraut} (Proposition XIV.4.1), and this is the only hypothesis needed for Godement's theorem.  

With this in mind, the proof is just an application of the $\mathrm{L}^1$--\,Godement theorem as stated in~\citep{kernelspaper}. Precisely, $f$ is positive definite only if $\hat{f}$ is positive and integrable (which here means it satisfies the integrability condition (\ref{eq:hc_integrable})). Then, $f$ is given by the inversion formula (\ref{eq:spherical_inverse}), which is the same as (\ref{eq:L1G}) with $g = C_{\scriptscriptstyle N}\hspace{0.02cm}\hat{f}$.

To finish the proof, note that uniqueness of $g$ follows by injectivity of the inverse spherical transform (the linear map that takes $\hat{f}$ to $f$ according to (\ref{eq:spherical_inverse})).

\section{Proof of Proposition \ref{prop:langtrick}} \label{app:langtrick}
\noindent \textbf{Part (a)\,:} in order to prove (\ref{eq:PD1}), it is enough to show that, on the right-hand side of (\ref{eq:L1G_bis}),
\begin{equation} \label{eq:proofPD1_1}
 \frac{1}{N!}\int_{\mathbb{R}^{\scriptscriptstyle N}}\,g(t)
V(t)\det\!\left[\rho^{\mathrm{i}t_{\ell}+\scriptscriptstyle (N-1)/2}_k\right]\hspace{0.02cm}dt = 
(\det(x))^{\scriptscriptstyle (N-1)/2}\,\det\!\left[\int_{\mathbb{R}}\,\gamma(t)t^{k-1}\hspace{0.03cm}e^{\mathrm{i}t s_\ell}\hspace{0.03cm}dt \right]
\end{equation}
where $g(t) = \gamma(t_{\scriptscriptstyle 1})\ldots\gamma(t_{\scriptscriptstyle N})$ and $s_\ell = \log(\rho_\ell)$. Note first that
$$
\det\!\left[\rho^{\mathrm{i}t_{\ell}+\scriptscriptstyle (N-1)/2}_k\right] = \prod^N_{k=1} \rho^{\scriptscriptstyle (N-1)/2}_k\times\det\!\left[\rho^{\mathrm{i}t_{\ell}}_k\right] = (\det(x))^{\scriptscriptstyle (N-1)/2}\times\det\!\left[\rho^{\mathrm{i}t_{\ell}}_k\right]
$$
This implies that (\ref{eq:proofPD1_1}) is equivalent to 
$$
 \frac{1}{N!}\int_{\mathbb{R}^{\scriptscriptstyle N}}\,g(t)
V(t)\det\!\left[\rho^{\mathrm{i}t_{\ell}}_k\right]\hspace{0.02cm}dt = 
\det\!\left[\int_{\mathbb{R}}\,\gamma(t)t^{k-1}\hspace{0.03cm}e^{\mathrm{i}t s_\ell}\hspace{0.03cm}dt \right]
$$
or, what is the same if $V(t)$ is expressed as a determinant,
\begin{equation} \label{eq:proofPD1_2}
 \frac{1}{N!}\int_{\mathbb{R}^{\scriptscriptstyle N}}\,g(t)
\det\!\left[t^{k-1}_{\ell}\right]\det\!\left[\rho^{\mathrm{i}t_{\ell}}_k\right]\hspace{0.02cm}dt = 
\det\!\left[\int_{\mathbb{R}}\,\gamma(t)t^{k-1}\hspace{0.03cm}e^{\mathrm{i}t s_\ell}\hspace{0.03cm}dt \right]
\end{equation}
Here, using the Andr\'eief identity~\citep{livan} (Chapter 11, Page 75), the left-hand side is equal to
$$
\frac{1}{N!}\int_{\mathbb{R}^{\scriptscriptstyle N}}\,
\det\!\left[t^{k-1}_{\ell}\right]\det\!\left[\rho^{\mathrm{i}t_{\ell}}_k\right]\hspace{0.02cm}\prod^N_{\ell=1} \gamma(t_\ell)dt_\ell = 
\det\!\left[\int_{\mathbb{R}}\,\gamma(t)\hspace{0.02cm}t^{k-1}\hspace{0.03cm}\rho^{\mathrm{i}t}_\ell\hspace{0.03cm}dt \right]
$$
which is the same as the right-hand side (by definition of $s_{\ell}$). Thus, (\ref{eq:proofPD1_2}) (equivalent to (\ref{eq:proofPD1_1})) has been proven true.

\noindent \textbf{Part (b)\,:} comparing (\ref{eq:L1G_bis}) and (\ref{eq:PD2}), it becomes clear that one must show
\begin{equation} \label{eq:proofPD2_1}
 \frac{1}{N!}\int_{\mathbb{R}^{\scriptscriptstyle N}}\,g(t)
V(t)\det\!\left[\rho^{\mathrm{i}t_{\ell}+\scriptscriptstyle (N-1)/2}_k\right]\hspace{0.02cm}dt =   
(\det(x))^{\scriptscriptstyle {(N-1)/2}}\;\mathrm{i}^{\scriptscriptstyle N(N-1)/2}\;V\left(-\partial/\partial s\right)\tilde{g}(s)
\end{equation}
However, as in the proof of Part (a), this is equivalent to
\begin{equation} \label{eq:proofPD2_2}
 \frac{1}{N!}\int_{\mathbb{R}^{\scriptscriptstyle N}}\,g(t)
V(t)\det\!\left[\rho^{\mathrm{i}t_{\ell}}_k\right]\hspace{0.02cm}dt =   
 \mathrm{i}^{\scriptscriptstyle N(N-1)/2}\;V\left(-\partial/\partial s\right)\tilde{g}(s)
\end{equation}
After writing the Leibniz expansion of the determinant, the left-hand side is equal to
$$
\frac{1}{N!}\sum_{w \in S_{\scriptscriptstyle N}} \varepsilon(w)
\int_{\mathbb{R}^{\scriptscriptstyle N}}\,g(t)
V(t)\prod^N_{k=1} \rho^{\mathrm{i}t_{w(k)}}_k
\hspace{0.02cm}dt =  \frac{1}{N!}\sum_{w \in S_{\scriptscriptstyle N}} \varepsilon(w)
\int_{\mathbb{R}^{\scriptscriptstyle N}}\,g(t)
V(t) e^{\mathrm{i}(s,w\hspace{0.02cm}t)}
\hspace{0.02cm}dt 
$$
Here, $S_{\scriptscriptstyle N}$ is the group of permutations of $\lbrace 1,\ldots, N\rbrace$ and $\varepsilon(w)$ is the signature of the permutation $w$. Moreover, on~the right-hand side $s_\ell = \log(\rho_\ell)$ and $w\hspace{0.02cm}t$ denotes the action of the permutation $w$ on $(t_{\scriptscriptstyle 1},\ldots, t_{\scriptscriptstyle N})$. By introducing a new variable of integration $u = w\hspace{0.02cm}t$ in each one of the integrals under the sum,
$$
\frac{1}{N!}\sum_{w \in S_{\scriptscriptstyle N}} \varepsilon(w)
\int_{\mathbb{R}^{\scriptscriptstyle N}}\,g(t)
V(t) e^{\mathrm{i}(s,w\hspace{0.02cm}t)}
\hspace{0.02cm}dt  = 
\frac{1}{N!}\sum_{w \in S_{\scriptscriptstyle N}} \varepsilon(w)
\int_{\mathbb{R}^{\scriptscriptstyle N}}\,g(w^{-\scriptscriptstyle 1}\hspace{0.02cm}u)
V(w^{-\scriptscriptstyle 1}\hspace{0.02cm}u) e^{\mathrm{i}(s,u)}
\hspace{0.02cm}du
$$
But the function $g$ is symmetric, while the Vandermonde polynomial $V$ is alternating --- for any permutation $w$, $g(w\hspace{0.02cm}u) = g(u)$ and $V(w\hspace{0.02cm}u) = \varepsilon(w)V(u)$. Therefore, the above sum is equal to
$$
\frac{1}{N!}\sum_{w \in S_{\scriptscriptstyle N}}
\int_{\mathbb{R}^{\scriptscriptstyle N}}\,g(u)
V(u) e^{\mathrm{i}(s,u)}
\hspace{0.02cm}du = \int_{\mathbb{R}^{\scriptscriptstyle N}}\,g(u)
V(u) e^{\mathrm{i}(s,u)}
\hspace{0.02cm}du
$$
and it now follows that the left-hand side of (\ref{eq:proofPD2_2}) is
\begin{equation} \label{eq:proofPD2_4}
 \frac{1}{N!}\int_{\mathbb{R}^{\scriptscriptstyle N}}\,g(t)
V(t)\det\!\left[\rho^{\mathrm{i}t_{\ell}}_k\right]\hspace{0.02cm}dt =   
\int_{\mathbb{R}^{\scriptscriptstyle N}}\,g(u)
V(u) e^{\mathrm{i}(s,u)}\hspace{0.02cm}du
\end{equation}
Finally, recalling the definition of the inverse Fourier transform $\tilde{g}(s)$, and differentiating under the integral, one has
$$
V(\partial/\partial s)\tilde{g}(s) =  \mathrm{i}^{\scriptscriptstyle N(N-1)/2}\;\int_{\mathbb{R}^{\scriptscriptstyle N}}\,g(u)
V(u) e^{\mathrm{i}(s,u)}\hspace{0.02cm}du
$$
which can be replaced back into (\ref{eq:proofPD2_4}) to obtain
$$
 \frac{1}{N!}\int_{\mathbb{R}^{\scriptscriptstyle N}}\,g(t)
V(t)\det\!\left[\rho^{\mathrm{i}t_{\ell}}_k\right]\hspace{0.02cm}dt =
(-\mathrm{i})^{\scriptscriptstyle N(N-1)/2}\;   
V(\partial/\partial s)\tilde{g}(s) = \mathrm{i}^{\scriptscriptstyle N(N-1)/2}\;   
V(-\partial/\partial s)\tilde{g}(s)
$$
which is identical to (\ref{eq:proofPD2_2}), as required. 

\section{Proof of Proposition \ref{prop:ramanujan}} \label{app:ramanujan}
The proof relies heavily on Ramanujan's theorem for symmetric cones, as given in~\citep{ding} (Page 450). This states that if the coefficients $a(m)$ in (\ref{eq:spherical_taylor}) are of the form $a(m) = q(m-\delta)$ for a function $q : \mathbb{C}^{\scriptscriptstyle N} \rightarrow \mathbb{C}$, which is symmetric, holomorphic for $\mathrm{Re}(\lambda_k) > -L$, where $L >3(N-1)/2$, and satisfies the growth condition
\begin{equation} \label{eq:growth_bis}
 \left|q(\lambda)\right| \leq C_{\scriptscriptstyle N}\hspace{0.02cm}\left|\Gamma_M(N + \lambda +\delta) \right|\prod^N_{k=1} e^{P\hspace{0.02cm}\mathrm{Re}(\lambda_k)}\times 
    e^{A\hspace{0.02cm}|\mathrm{Im}(\lambda_k)|} 
\end{equation}
where $P,A > 0$ and $A <\pi$, then the following hold.
\begin{itemize}
    \item[(a)] the series (\ref{eq:spherical_taylor}) converges in a neighborhood of $x = 0$, where it defines a real-analytic function $F$.
    \item[(b)] this function extends continuously to all of $M$, by the following absolutely convergent integral
    \begin{equation} \label{eq:ding1}
        F(x) = C_{\scriptscriptstyle N}\,\int_{\mathbb{R}^{\scriptscriptstyle N}} \tilde{F}(\sigma+\mathrm{i}t)\hspace{0.03cm}\Phi_{\sigma +\mathrm{i}t}(x)\hspace{0.03cm}(V(t))^2dt 
    \end{equation}
where $-L + (N-1)/2 <\sigma < -(N-1)$ and $\tilde{F}(\lambda) = \Gamma_M(\delta-\lambda)\hspace{0.02cm}q(\lambda)$.
\item[(c)] for any $\sigma$ as above, $\omega$ in the convex hull of $w\hspace{0.02cm}\delta$, where $w$ ranges over the symmetric group $S_{\scriptscriptstyle N\hspace{0.03cm}}$,
\begin{equation} \label{eq:ding2}
  \tilde{F}(\lambda) = \int_M F(x)\hspace{0.03cm}\Phi_{-\lambda}(x)\hspace{0.03cm}\mathrm{vol}(dx)     
\end{equation}
is an absolutely convergent integral, whenever $\lambda = \sigma + \omega + \mathrm{i}t$ with $t \in \mathbb{R}^{\scriptscriptstyle N}$. 
\end{itemize}

Here, the notation $\sigma + \lambda$, where $\sigma$ is a number and $\lambda \in \mathbb{C}^{\scriptscriptstyle N}$, means $\sigma$ is added to each component of $\lambda$. 

For the proof of Proposition \ref{prop:ramanujan}, in view of (\ref{eq:rama_cond1}), take
\begin{equation} \label{eq:proofQQ}
    q(\lambda) = \Gamma_M(2\alpha+\lambda+\delta)\hspace{0.02cm}\psi(\lambda)
\end{equation}
which is then symmetric, because $\psi$ is symmetric by assumption and $\Gamma_M(\lambda + \delta)$ is a symmetric function of $\lambda$, as can be seen after replacing $\delta_k = k -(N+1)/2$ into (\ref{eq:GammaM}), and holomorphic for $\mathrm{Re}(\lambda_k) > -L$, where $L = 2\alpha - (N-1)/2 > 3(N-1)/2$ (recall that $\alpha > N-1$), again by the assumptions made regarding $\psi$ and the fact that the Gamma function $\Gamma_M(\lambda)$ is holomorphic for $\mathrm{Re}(\lambda_k) > N-1$, as one may see from (\ref{eq:GammaM}). Moreover, this $q$ satisfies the growth condition (\ref{eq:growth_bis}), since (\ref{eq:rama_cond2}) implies
\begin{equation} \label{eq:growth_proof_1}
|q(\lambda)| 
\leq C_{\scriptscriptstyle N}\hspace{0.02cm}\left|\Gamma_M(2\alpha + \lambda +\delta) \right|\prod^N_{k=1} e^{P\hspace{0.02cm}\mathrm{Re}(\lambda_k)}\times 
    e^{A\hspace{0.02cm}|\mathrm{Im}(\lambda_k)|}
\end{equation}
and since, from (\ref{eq:GammaM}), after putting $z_k = \lambda_k + 2\alpha - (N-1)/2$,
\begin{equation} \label{eq:growthèproof_2}
\frac{\Gamma_M(2\alpha+\lambda+\delta)}{\Gamma_M(N+\lambda+\delta)}
= \prod^N_{k=1} \frac{\Gamma(z_k)}{\Gamma(z_k + N - 2\alpha)} = \prod^N_{k=1} z^{2\alpha - N}_k\left(1+O(|z_k|^{-1})\right)
\end{equation}
where the second equality follows from the asymptotic form for the ratio of two Gamma functions~\citep{olver} (Page 119). Indeed, by (\ref{eq:growth_proof_1}) and (\ref{eq:growthèproof_2}),
$$
|q(\lambda)| 
\leq C_{\scriptscriptstyle N}\hspace{0.02cm}\times \prod^N_{k=1} z^{2\alpha - N}_k\left(1+O(|z|^{-1})\right)\times \left|\Gamma_M(N + \lambda +\delta) \right|\prod^N_{k=1} e^{P\hspace{0.02cm}\mathrm{Re}(\lambda_k)}\times 
    e^{A\hspace{0.02cm}|\mathrm{Im}(\lambda_k)|}
$$
and this implies (\ref{eq:growth_bis}) because $|z|^{\hspace{0.02cm}2\alpha - N} = o\left(e^{\hspace{0.02cm}\varepsilon\hspace{0.02cm}(\mathrm{Re}(z)+\mathrm{Im}(z))}\right)$ for any $\varepsilon > 0$, in the limit where $|z|\rightarrow \infty$. 
Thus, the conditions of Ramanujan's master theorem are all verified for $q$ given by (\ref{eq:proofQQ}), so items (a) to (c) above can be applied in the context of Proposition \ref{prop:ramanujan}. First, (a) implies the spherical series (\ref{eq:spherical_taylor}) defines a real-analytic function $F$ in the neighborhood of $x = 0$. Moreover, (b) and (c) imply that this $F$ extends continuously to all of $M$, in such a way that (\ref{eq:ding1}) and (\ref{eq:ding2}) are satisfied. Let $\sigma = -\alpha$ in (\ref{eq:ding1}) and note that (\ref{eq:sphericalfunction}) implies $\Phi_{\mathrm{i}t-\alpha}(x) = \Delta^{-\alpha}(x)\Phi_{\mathrm{i}t}(x)$. Then, (\ref{eq:ding1}) and the subsequent definition of $\tilde{F}(\lambda)$, applied with $q$ as in (\ref{eq:proofQQ}), show that
\begin{equation} \label{eq:ding11}
  \Delta^\alpha(x)F(x) = 
  C_{\scriptscriptstyle N}\,\int_{\mathbb{R}^{\scriptscriptstyle N}} \left|\Gamma_M(\alpha + \delta + \mathrm{i}t)\right|^2\psi(\mathrm{i}t - \alpha)\hspace{0.03cm}\Phi_{\mathrm{i}t}(x)\hspace{0.03cm}(V(t))^2dt 
\end{equation}
Now, letting $f(x) = \Delta^\alpha(x)F(x)$, it becomes clear $f$ is $U$-invariant, since each spherical function $\Phi_{\mathrm{i}t}$ is $U$-invariant. On the other hand, choosing $\sigma = -\alpha$, $\omega = \delta$, and $t = 0$ in (\ref{eq:ding2}), and noting that $\Phi_{-\delta}(x) = 1$ is a constant function, it follows that $f$ is integrble. The spherical transform of $f$ is found from (\ref{eq:ding2}), with $\lambda = -\alpha+\mathrm{i}t$, which shows that this spherical transform is equal to $\hat{f}$ in (\ref{eq:ramanujan}). 

Before pursuing the final stage of the proof, note that $\hat{f}$ (given by (\ref{eq:ramanujan})) satisfies the integrability condition (\ref{eq:hc_integrable}). Indeed, (\ref{eq:rama_cond2}) and (\ref{eq:ramanujan}) imply that
$$
|\hat{f}(t)| = 
\left|\Gamma_M(\alpha + \delta + \mathrm{i}t)\right|^2\left|\psi(\mathrm{i}t - \alpha)\right| \leq
\left|\Gamma_M(\alpha + \delta + \mathrm{i}t)\right|^2 \times
C_{\scriptscriptstyle N}e^{-NP\hspace{0.02cm}\alpha}\hspace{0.02cm}\prod^N_{k=1}  
    e^{A\hspace{0.02cm}|t_k|}
$$
However, by (\ref{eq:GammaM}) and~\citep{beals} (Corollary 2.5.3)
\begin{flalign}
\nonumber \left|\Gamma_M(\alpha + \delta + \mathrm{i}t)\right|^2 & = (2\pi)^{N(N-1)}\prod^N_{k=1} \left| \Gamma\left(\mathrm{i}t_k + \alpha - (N-1)/2\right)\right|^{\hspace{0.02cm} 2} \\[0.1cm] 
\nonumber & =
(2\pi)^{N^2}\prod^N_{k=1}  |t_k|^{\alpha - N/2}\hspace{0.03cm}e^{-\pi|t_k|}\left(1+O(|t_k|^{-1})\right)    
\end{flalign}
and it therefore follows 
\begin{equation} \label{eq:AandPI}
|\hat{f}(t)| \leq 
C_{\scriptscriptstyle N}\,e^{-NP\hspace{0.02cm}\alpha}\hspace{0.02cm}\prod^N_{k=1} |t_k|^{\alpha - N/2}\hspace{0.03cm}e^{(A-\pi)|t_k|}
\end{equation}
But since $A < \pi$, this is exponentially small as the $|t_k|$ become large. In turn, this shows that (\ref{eq:hc_integrable}) is satisfied. 

To complete the proof of Proposition \ref{prop:ramanujan}, it remains to show that $F$ in (\ref{eq:ding11}) is real-analytic (so that it gives an analytic extension, rather than just the continuous extension stated in item (b)), and also that $f$ is positive definite if and only if $\psi(\mathrm{i}t - \alpha) \geq 0$ for all $t \in \mathbb{R}^{\scriptscriptstyle N}$. To do so, rewrite (\ref{eq:ding11}) as follows
$$
\Delta^\alpha(x)F(x) = C_{\scriptscriptstyle N}\,\int_{\mathbb{R}^{\scriptscriptstyle N}} \hat{f}(t)\hspace{0.03cm}\Phi_{\mathrm{i}t}(x)\hspace{0.03cm}(V(t))^2dt = 
\frac{C_{\scriptscriptstyle N}}{V(\mathrm{i}\rho)}\,\int_{\mathbb{R}^{\scriptscriptstyle N}}\hat{f}(t)
V(t)\det\!\left[\rho^{\mathrm{i}t_{\ell}+\scriptscriptstyle (N-1)/2}_k\right]\hspace{0.02cm}dt
$$
where the second equality follows from (\ref{eq:spherical_inverse}) and (\ref{eq:spherical_inverse_bis}) (here, the factor $1/N!$ from (\ref{eq:spherical_inverse_bis}) has been absorbed into $C_{\scriptscriptstyle N}$). This easily simplifies to
\begin{equation} \label{eq:holom1}
\Delta^{\alpha-(N-1)/2}(x)\hspace{0.03cm}F(x) =
\frac{C_{\scriptscriptstyle N}}{V(\mathrm{i}\rho)}\,\int_{\mathbb{R}^{\scriptscriptstyle N}}\,\hat{f}(t)
V(t)\det\!\left[\rho^{\mathrm{i}t_{\ell}}_k\right]\hspace{0.02cm}dt
\end{equation}
Now, replace $\rho_k$ with the complex variable $\chi_{k} = \rho_k\hspace{0.03cm}e^{\mathrm{i}\hspace{0.02cm}\varphi_k}$ where $\varphi_k$ is real and $|\varphi_k| <  \varepsilon$ with $\varepsilon <$ the minimum~of $\pi/2$ and $(\pi-A)/2$ (here, $A < \pi$ is the constant in (\ref{eq:AandPI})). With this replacement, the right-hand side of (\ref{eq:holom1}) becomes
\begin{equation} \label{eq:holom3}
\Psi(\chi) = \frac{C_{\scriptscriptstyle N}}{V(\mathrm{i}\chi)}\,\int_{\mathbb{R}^{\scriptscriptstyle N}}\,\hat{f}(t)
V(t)\det\!\left[\chi^{\mathrm{i}t_{\ell}}_k\right]\hspace{0.02cm}dt 
\end{equation}
The aim is to show that $\Psi$ is a holomorphic function of $\chi$ (each $\chi_k$ having its argument $\varphi_k$ subject to $|\varphi_k| < \varepsilon$). It will then follow that (\ref{eq:holom1}) defines a real-analytic function of $\rho$ (the restriction of a holomorphic function of $\chi$). Moreover, since $\Delta^{\alpha-(N-1)/2}(x)$ is analytic and non-zero, as a function of $\rho$, this will show that $F$ is real-analytic. 
To see that $\Psi$ is holomorphic, note that the integral in (\ref{eq:holom3}) defines a holomorphic function of $\chi$ -- call this $I(\chi)$. Indeed, the~function under that integral is holomorphic for each fixed $t$, and the integral converges uniformly in $\chi$, due to the following upper bound (the fact that the integral in (\ref{eq:holom3}) is holomorphic then follows from~\citep{bochner} (Page 41))
$$
\left|\hspace{0.03cm}\hat{f}(t)
V(t)\det\!\left[\chi^{\mathrm{i}t_{\ell}}_k\right]\right| \leq 
|P(t)|\hspace{0.02cm}\prod^N_{k=1} e^{(A-\pi)|t_k|}\times \left| \det\!\left[\chi^{\mathrm{i}t_{\ell}}_k\right] \right| \leq 
|P(t)|\hspace{0.02cm}\prod^N_{k=1} e^{(A-\pi+\varepsilon)|t_k|}
$$
where $P(t)$ is some polynomial and $A-\pi+\varepsilon < 0$. Here, the first inequality follows from (\ref{eq:AandPI}), and the second one since $|\chi^{\mathrm{i}t_{\ell}}_k| = e^{-t_{\ell}\hspace{0.02cm}\varphi_k}$. Therefore, $\Psi(\chi) = C_{\scriptscriptstyle N}\hspace{0.02cm}I(\chi)/V(\mathrm{i}\chi)$ is a ratio of holomorphic functions, and is holomorphic because any zero of $V(\mathrm{i}\chi)$ is also a zero of $I(\chi)$ (note that zeros of  $V(\mathrm{i}\chi)$ occur when $\chi_k = \chi_\ell$ for some $k < \ell$).

The final part of the proof requires showing that $f$ (defined after (\ref{eq:ding11})) is positive definite if and only if  
 $\psi(\mathrm{i}t-\alpha) \geq 0$ for all $t \in \mathbb{R^{\scriptscriptstyle N}}$. Recall that $f$ is integrable and its spherical transform is $\hat{f}$ given in (\ref{eq:ramanujan}). Recall also that $\hat{f}$ satisfies the integrability condition (\ref{eq:hc_integrable}), as shown above, using (\ref{eq:AandPI}). Therefore, $f$ is given by the inversion formula (\ref{eq:spherical_inverse}). Accordingly, from (\ref{eq:spherical_inverse}) and (\ref{eq:L1G}) of Theorem \ref{th:L1G}, by injectivity of the inverse spherical transform, $f$ is positive definite if and only if $\hat{f}(t) \geq 0$ for all $t \in \mathbb{R}^{\scriptscriptstyle N}$. This happens precisely when $\psi(\mathrm{i}t-\alpha) \geq 0$ for all $t \in \mathbb{R^{\scriptscriptstyle N}}$ (this is by (\ref{eq:ramanujan})).   
 
 




\end{document}